\pdfoutput=1
\RequirePackage{ifpdf}
\ifpdf 
\documentclass[pdftex]{sigma}
\else
\documentclass{sigma}
\fi

\newtheorem{thm}{Theorem}[section]

\newtheorem{lem}[thm]{Lemma}
\newtheorem{prop}[thm]{Proposition}
{\theoremstyle{definition} \newtheorem{defn}[thm]{Definition}
\newtheorem{rem}[thm]{Remark}  }

\numberwithin{equation}{section}
\DeclareMathOperator{\id}{id}

\begin{document}

\allowdisplaybreaks

\renewcommand{\PaperNumber}{001}

\FirstPageHeading

\ShortArticleName{On Classif\/ication of Finite-Dimensional Superbialgebras and Hopf Superalgebras}

\ArticleName{On Classif\/ication of Finite-Dimensional\\
Superbialgebras and Hopf Superalgebras}

\Author{Said AISSAOUI~$^\dag$ and Abdenacer MAKHLOUF~$^\ddag$}

\AuthorNameForHeading{S.~Aissaoui and A.~Makhlouf}

\Address{$^\dag$~Universit\'{e} A-Mira, Laboratoire de Math\'{e}matiques Appliqu\'{e}es,\\
\hphantom{$^\dag$}~Targa Ouzemmour 06000 B\'{e}jaia, Algeria}
\EmailD{\href{mailto:aissaoui.said@yahoo.fr}{aissaoui.said@yahoo.fr}}

\Address{$^\ddag$~Universit\'{e} de Haute Alsace, Laboratoire de Math\'{e}matiques, Informatique et Applications,\\
\hphantom{$^\ddag$}~4, rue des Fr\`{e}res Lumi\`{e}re F-68093 Mulhouse, France}
\EmailD{\href{mailto:Abdenacer.Makhlouf@uha.fr}{Abdenacer.Makhlouf@uha.fr}}
\URLaddressD{\url{http://www.algebre.fst.uha.fr/makhlouf/Makhlouf.html}}

\ArticleDates{Received February 08, 2013, in f\/inal form December 23, 2013; Published online January 02, 2014}

\Abstract{The purpose of this paper is to investigate f\/inite-dimensional superbialgebras and Hopf
superalgebras.
We study connected superbialgebras and provide a~classif\/ication of non-trivial superbialgebras and Hopf
superalgebras in dimension $n$ with $n\leq 4$.}


\Classification{16T10; 57T05; 17A70}

\pdfbookmark[1]{Introduction}{intro}
\section*{Introduction}

Hopf superalgebras are a~generalization of the supergroup notion.
In physics, inf\/inite-di\-men\-sio\-nal Hopf superalgebras based on Lie superalgebras turn out to be related to
the integrable $S$-matrix of the AdS/CFT correspondence.
For a~large class of f\/inite-dimensional Lie superalgebras (including the classical simple ones) a~Lie
supergroup associated to the algebra is def\/ined by f\/ixing the Hopf superalgebra of functions on the
supergroup~\cite{sheunert}.
This framework seems to be interesting for studying supermanifolds and supersymmetry.
Majid showed that the theo\-ry of Lie superalgebras and Hopf superalgebras can be reduced to the classical
case using the bosonization by $\mathbb{K}[\mathbb{Z}/2\mathbb{Z}]$, where $\mathbb{K}$ is an algebraically
closed f\/ield~\cite{Majid2}.
Recently some classes of Hopf superalgebras were investigated like for example pointed Hopf superalgebras
and quasitriangular Hopf superalgebras~\cite{Andrus-Schneider2,Andrus2,GOULD}, see also~\cite{deguchi1990,kulish1989}.
So far very few is known about general classif\/ication of superbialgebras and Hopf superalgebras.
Nevertheless, classif\/ication of f\/inite-dimensional Hopf algebras is known for small dimensions and
structure of many classes are deeply studied, see for
example~\cite{Andrus1,Andrus2,Andrus-Schneider3,Beattie,Beattie3,Etingof-Gelaki,fukuda,Masuoka,Milnor,Montgomery,natale,Ng,Stefan2,stefan,williams,Zhu}.
For general theory of Hopf algebras, we refer to~\cite{Ab,Guichardet,Kassel,Majid,MontgomeryLivre,Shnider},
and for applications see for example~\cite{CK98,Holtkamp,Kreimer97,Makhlouf-Hermann}.

In this paper, we discuss properties of $n$-dimensional superbialgebras and provide a~classif\/ication of
non-trivial superbialgebras in dimensions $2$, $3$ and $4$.
Moreover we derive a~classif\/ication of Hopf superalgebras for these dimensions.

The paper is organized as follows.
Section~\ref{review} reviews def\/initions and properties of superbialgebras and Hopf superalgebras.
Moreover, we provide an outline of the computations leading to classif\/ications.
In Section~\ref{trivconnsupalg}, we focus on trivial and connected superbialgebras.
We cha\-racterize connected superbialgebras and classify connected $2$-dimensional superbialgebras and Hopf
superalgebras.
Sections~\ref{classification3} and~\ref{classification4} deal respectively
with classif\/ication of $3$- and $4$-dimensional superbialgebras and Hopf superalgebras.
We emphasize on non-trivial superbialgebras, from which we derive the non-trivial Hopf superalgebras.
To this end, we establish classif\/ication of $3$-dimensional superalgebras and use the $4$-dimensional
classif\/ication due to Armour, Chen and Zhang~\cite{Armour1}.

Throughout this paper, we work over $\mathbb K$, an algebraically closed f\/ield of characteristic zero.
Unless otherwise specif\/ied, we will denote the multiplication by $\mu$ or $\cdot$, and the unit element
by $\eta(1)=e_{1}=e_{1}^{0}=1$.
For simplicity, we will use just product (or multiplication), unit, counit and coproduct (or
comultiplication) in the sense of even morphism (i.e.\ we will drop the suf\/f\/ix `super').
We also would like to mention that in the last section, most of the results are obtained by using the
computer algebra system Mathematica.
The program we have used is available.

\section{Definitions of superbialgebras and Hopf superalgebras }
\label{review}

In this section, we summarize def\/initions and properties of superbialgebras and Hopf superalgebras.
For more details, we refer to~\cite{Andrus-Schneider2,GOULD}.

A \textit{superspace} $A$ is a~$\mathbb{K}$-vector space endowed with a~$\mathbb{Z}/$2$\mathbb{Z}$-grading,
in other words, it writes as a~direct sum of two vector spaces $A=A_{0} \oplus A_{1}$ such as $A_{0}$ is
the even part and $A_{1}$ is the odd part.
Elements of $A_{0}$ (resp.\ $A_{1}$) are called \emph{even homogeneous} (resp.\ \emph{odd homogeneous}).
If $a \in A_{0}$, we set $|a|=\deg(a)=0$ and if $a\in A_{1}$, $|a|=\deg(a)=1$.
Notice that some authors denote the dimension of $A$ by $i|j$, such that $i=n_{0}$ and $j=n_{1}$, where
$n_{0}=\dim A_{0}$ and $n_{1}=\dim A_{1}$.
But in this paper we preserve this notation for the superalgebra obtained by the $i$th algebra.
Notice that the following def\/initions could be described in a~symmetric category of superspaces.
\begin{defn}
A \textit{superalgebra} is a~triple $(A,\mu,\eta)$ where $A$ is a~superspace, $\mu: A \otimes A \rightarrow A$
(multiplication) and $\eta:\mathbb K\rightarrow A$ (unit) are two superspace morphisms satisfying
\begin{gather}
\label{asso10}
\mu\circ(\mu\otimes \id_{A})=\mu\circ(\id_{A}\otimes\mu)
\qquad
\text{(associativity)},
\\
\label{unit10}
\mu\circ\eta\otimes \id_{A}=\mu\circ \id_{A}\otimes\eta
\qquad
\text{(unity)}.
\end{gather}
\end{defn}

A superalgebra $A$ is called \textit{commutative} if the multiplication satisf\/ies $\mu \circ \tau=\mu $
where $\tau$ is a~superf\/lip.
In other words, for all homogeneous elements $a, b \in A$, $\mu(a \otimes b)=(-1)^{|a||b|}\mu(b \otimes a)$.

Let $(A, \mu_{A}, \eta_{A})$, $(B, \mu_{B}, \eta_{B})$ be two superalgebras.
A map $f:A\rightarrow B$ is a~\textit{superalgebra morphism} if it is a~superspace morphism satisfying
\begin{gather}
\label{MorphAlg}
f\circ\mu_{A}=\mu_{B}\circ f\otimes f
\qquad
\text{and}
\qquad
f\circ\eta_{A}=\eta_{B}.
\end{gather}
\begin{defn}
A \textit{supercoalgebra} is a~triple $(C,\Delta,\varepsilon)$
where $C$ is a~superspace such that $\Delta: C \rightarrow C \otimes C$ (comultiplication or coproduct)
and $\varepsilon: C \rightarrow\mathbb K$ (counit) are two superspace morphisms satisfying
\begin{gather}
\label{coasso10}
(\Delta\otimes \id_C)\circ\Delta=(\id_C\otimes\Delta)\circ\Delta
\qquad
\text{(coassociativity),}
\\
\label{coasso20}
(\varepsilon\otimes \id_{C})\circ\Delta=(\id_{C}\otimes\varepsilon)\circ\Delta
\qquad
\text{(counity).
}
\end{gather}

\end{defn}
We use the Sweedleer's notation for the coproduct, we set for all $x \in C $\begin{gather*}
\Delta(x)=\sum_{(x)}x^{(1)}\otimes x^{(2)},
\qquad
(\id\otimes\Delta)\circ\Delta(x)=(\Delta\otimes \id)\circ\Delta(x)=\sum_{(x)}x^{(1)}\otimes x^{(2)}
\otimes x^{(3)}.
\end{gather*}

A supercoalgebra $C$ is said to be \textit{cocommutative} if the comultiplication satisf\/ies $\Delta=\tau
\circ \Delta$, where $\tau$ is the superf\/lip, that is
\begin{gather*}
\forall\, x\in C,
\qquad
\Delta(x)=\sum_{(x)}x^{(1)}\otimes x^{(2)}=\sum_{(x)}(-1)^{|x^{(1)}||x^{(2)}|}x^{(2)}\otimes x^{(1)}.
\end{gather*}

Let $(A, \Delta_{A}, \varepsilon_{A})$, $(B,\Delta_{B}, \varepsilon_{B})$ be two supercoalgebras.
The map $f:(A, \Delta_{A}, \varepsilon_{A})\rightarrow (B, \Delta_{B}, \varepsilon_{B})$ is
a~\textit{supercoalgebra morphism} if it is a~superspace morphism and satisf\/ies
\begin{gather}
\label{morp}
f\otimes f\circ\Delta_{A}=\Delta_{B}\circ f
\qquad
\text{and}
\qquad
\varepsilon_{B}\circ f=\varepsilon_{A}.
\end{gather}

Now, we consider a~superbialgebra structure which is obtained by combining superalgebras and
supercoalgebras.
Moreover, if a~superbialgebra has an anti-superbialgebra morphism satisfying some condition, it is called
a~Hopf superalgebra.
\begin{defn}
A \textit{superbialgebra} is a~tuple $(A, \mu, \eta, \Delta, \varepsilon)$, where $(A,\mu,\eta)$ is
a~superalgebra and $(A, \Delta, \varepsilon)$ is a~supercoalgebra such that one of these two equivalent
compatibility conditions hold:
\begin{enumerate}\itemsep=0pt
\item[1)] $\Delta: A \rightarrow A \otimes A $ and $\varepsilon: A \rightarrow \mathbb K$ are
superalgebra morphisms,
\item[2)] $\mu: A \otimes A \rightarrow A $ and $\eta:\mathbb K\rightarrow A$ are supercoalgebra morphisms.
\end{enumerate}
\end{defn}
In other words, $\Delta$ (resp.\ $\varepsilon $) satisf\/ies the following compatibility condition
\begin{gather}
\label{cop}
\Delta\circ\mu=(\mu\otimes\mu)\circ(\id_{A}\otimes\tau\otimes \id_{A})\circ(\Delta\otimes\Delta),
\qquad
\Delta\circ\eta=\eta\otimes\eta
\\
\nonumber
(\text{resp.}
\quad
\varepsilon\circ\mu=\mu_{\mathbb K}
\circ(\varepsilon\otimes\varepsilon),
\quad
\varepsilon\circ\eta=\id_{\mathbb K}),
\end{gather}
where $\mu_{\mathbb K}$ is the multiplication of $\mathbb K$.

The last condition in~\eqref{cop} says that the unit $1=e_{1}^{0}=\eta(1)$ of the superbialgebra is
a~group-like element $\Delta(1)=1 \otimes 1$.
\begin{rem}
Note that if $(A=A_{0}\oplus A_{1}, \mu, \eta, \Delta, \varepsilon)$ is a~superbialgebra with $\eta(1)=1$,
a~unit element, then we have
\begin{gather*}
\varepsilon(A_{1})={0}Ê
\qquad
\text{and}
\qquad
\varepsilon(1)=1.
\end{gather*}
Indeed, if $\varepsilon: A \longrightarrow \mathbb K$ is a~superspace morphism and $\mathbb K$ is
a~superspace whose odd part is {0}, then $\varepsilon$ sends the odd part of A to odd part of $\mathbb K$
which is ${0}$.
For the second assertion, if $A$ is a~superbialgebra then $\varepsilon$ is a~superalgebra morphism and it
sends the unit element of $A$ to unit element of $\mathbb K$.
\end{rem}
If $A$ and $B$ are two superbialgebras over $\mathbb K$, we shall call a~superspace morphism $ f: A
\rightarrow B $ a~\textit{superbialgebra morphism} if it is both a~superalgebra morphism and
a~supercoalgebra morphism.
\begin{defn}
A \textit{Hopf superalgebra} is a~superbialgebra admitting an antipode,
that is a~superspace morphism $S:A \rightarrow A$ which satisfies the following condition:
\begin{gather}
\label{antip}
\mu\circ(S\otimes \id_{A})\circ\Delta=\mu\circ(\id_{A}\otimes S)\circ\Delta=\eta\circ\varepsilon.
\end{gather}
A Hopf superalgebra is given by a~tuple $H=(A,\mu,\eta,\Delta,\varepsilon, S)$.
\end{defn}

Let $(A,\mu,\eta,\Delta,\varepsilon, S)$ be a~Hopf superalgebra and $S$ its antipode.
Then we have the following properties:
\begin{enumerate}\itemsep=0pt
\item[1)] $ S\circ \mu=\mu \circ \tau \circ (S\otimes S)$;
\item[2)] $S\circ \eta=\eta $;
\item[3)] $\varepsilon \circ S=\varepsilon $;
\item[4)] $\tau \circ (S\otimes S)\circ \Delta=\Delta\circ S $;
\item[5)] if
$A$ is commutative or cocommutative then $S \circ S=\id$ where $\id: A\longrightarrow A$ is the identity
morphism;

\item[6)] if the antipode $S$ of the Hopf superalgebra $A$, with $\Delta$ as comultiplication, is bijective
then $A$ is also a~Hopf superalgebra, with the opposite comultiplication $\Delta^{'}=\tau \circ \Delta$ and
the opposite antipode $S^{'}=S^{-1}$;
\item[7)] all f\/inite-dimensional Hopf superalgebras possess bijective
antipodes.
\end{enumerate}
\begin{rem}
Some superalgebras do not carry a~superbialgebra structure.
For example $A=M_{2}(\mathbb K)$~is a~superalgebra with matrix multiplication and such that the even part
is $A_{0}=\mathbb K e_{1}^{0}~\oplus~\mathbb K~e_{2}^{0}$ and the odd part is $A_{1}=\mathbb
K~e_{1}^{1}\oplus~\mathbb K~e_{2}^{1}$, where
\begin{gather*}
e_{1}^{0}=\left(
\begin{matrix}
1&0
\\
0&1
\\
\end{matrix}
\right),
\qquad
 e_{2}^{0}=\left(
\begin{matrix}
1&0
\\
0&0
\\
\end{matrix}
\right),
\qquad
 e_{1}^{1}=\left(
\begin{matrix}
0&1
\\
0&0
\\
\end{matrix}
\right),
\qquad
 e_{2}^{1}=\left(
\begin{matrix}
0&0
\\
1&0
\\
\end{matrix}
\right).
\end{gather*}

\end{rem}
Indeed, the compatibility condition between the counit $\varepsilon$ and the product
$\varepsilon(\mu(x\otimes y))=\varepsilon(x)\varepsilon(y)$ isn't always satisf\/ied, since in one hand we
have $\varepsilon (\mu(e_{1}^{1}\otimes e_{2}^{1}))=\varepsilon(e_{2}^{0})=0$, because
$ \varepsilon(e_{2}^{1})=\varepsilon(e_{1}^{1})=0$, and in the other hand
$\varepsilon(\mu(e_{2}^{1}\otimes e_{1}^{1}))=\varepsilon(e_{1}^{0}-e_{2}^{0})=0$.
Then $\varepsilon(e_{2}^{0})=1$, which leads to a~contradiction.

In Sections~\ref{classification3} and~\ref{classification4}, we aim to classify $n$-dimensional non-trivial
superbialgebras (resp.\ Hopf superalgebras), for $n=3$ and $n=4$.
An $n$-dimensional superbialgebra (resp.\ a Hopf superalgebra) is identif\/ied
to its structure constants with respect to a~f\/ixed basis.
It turns out that the axioms of superbialgebra structure translate to a~system of polynomial equations that
def\/ine the algebraic variety of $n$-dimensional superbialgebras which is embedded into $\mathbb
K^{2n_{0}^{3}+6n_{0}n_{1}^{2}+n_{0}-1}$.
The classif\/ication requires to solve this algebraic system.
The calculations are handled using a~computer algebra system.
We include in the following an outline of the computation.
\begin{enumerate}\itemsep=0pt
\item We provide a~list of non-trivial $n$-dimensional superalgebras:
\begin{itemize}\itemsep=0pt
\item For $n=2$, the list is taken from~\cite{D-Makh,Gabriel}.
\item For $n=3$, the list is given in Propositions~\ref{PropClassDim3-1} and~\ref{PropClassDim3-2}.
\item For $n=4$, the list is taken from~\cite{Armour}.
\end{itemize}
\item Fix a~superalgebra $A$ and compute the supercoalgebra's structure constants $\{D_{(i,s),(j,t)}^k\}$
and $\{\xi_{i}^{0}\}$, where
\begin{gather*}
\Delta (e_{i}^{l})=\sum_{s=0}^{1}\sum_{j=1}^{n_{s}}\sum_{k=1}^{n_{t}}{D_{(i, l)}^{(j,s)(k)}e_{j}^{s}\otimes e_{k}^{t}}
\end{gather*}
with $t=(l+s)\mod[2]$ and
\begin{gather*}
\varepsilon (e_{i}^{l})=
\begin{cases}
\xi_{i}^{0} & \text{if}
\quad
l=0
\quad
\text{and}
\quad
i\neq 1,
\\
0 & \text{if}
\quad
l=1,
\\
1 & \text{if}
\quad
l=0
\quad
\text{and}
\quad
i=1
\end{cases}
\end{gather*}
with respect to bases $\{e_{i}^{0}\}_{i=1,\dots, n_0}$ and $\{e_{i}^{1}\}_{i=1,\dots, n_1}$ of $A_0$
and $A_1$, respectively.
Supercoalgebra structures on $A$ making it a~superbialgebra are in one-to-one correspondence with solutions
of the system corresponding to~\eqref{coasso10}, \eqref{coasso20}, \eqref{cop}.
\item Obtained in (2) is a~family $\{A_i\}_{i\in I}$ of superalgebras whose underlying superalgebra is~$A$.
We sort them into isomorphism classes.
Two superbialgebras, given by their structure constants, are isomorphic if there exist matrices $(T_{(i,
0)}^{k},T_{(i, 1)}^{k})_{i,k}$ def\/ining a~superbialgebra morphism with respect to the basis, that is
satisfying $\det(T_{(i, 0)}^{k})\det(T_{(i, 1)}^{k})\neq 0$ and~\mbox{\eqref{MorphAlg}--\eqref{morp}}.

\item For each representative of isomorphism classes, determine whether it admits an antipode or not.
It follows from existence or not of structure constants $\lambda_{(i, s)}^{k}$, of a~morphism $S$,
$S(e_{i}^{s})=\sum\limits_{k=1}^{n_{s}}\lambda_{(i, s)}^{k}e_{k}^{s}$, satisfying~\eqref{antip}.
\end{enumerate}

\section{Trivial and connected superbialgebras}\label{trivconnsupalg}

We describe in this section some properties of connected superbialgebras, which are superbialgebras with
$1$-dimensional even part and trivial superbialgebras, which are superbialgebras with trivial odd part.
\subsection {Trivial superbialgebras} We have the following obvious result.
\begin{prop}
Every finite-dimensional bialgebra $($resp.\ Hopf algebra$)$ is a~superbialgebra $($resp.\
Hopf superalgebra$)$, whose odd part is reduced to $\{0\}$.
This superbialgebra $($resp.\ Hopf superalgebra$)$ is called trivial superbialgebra $($resp.\
trivial Hopf superalgebra$)$.
\end{prop}

The algebraic classif\/ication of Hopf algebra were investigated by many authors.
Hopf algebras of prime dimensions were classif\/ied in~\cite{Zhu}.
We refer to~\cite{stefan,williams} for the classif\/ication of f\/inite-dimensional Hopf algebras up to
dimension $11$, see also~\cite{Makhlouf-Hopf}.
Semisimple Hopf algebras of dimension $12$ were given in~\cite{fukuda} and complete classif\/ication of
dimension $12$ provided in~\cite{natale}.
See~\cite{Beattie2} for dimension $14$ and~\cite{Cheng-Ng} for dimension~$20$.
The minimal dimension where the problem is unsolved is~$24$.
For the classif\/ication of $2$ and $3$-dimensional bialgebras, we refer to~\cite{D-Makh}.
They are based on Gabriel's result of $2$ and $3$-dimensional algebras~\cite{Gabriel}.
Classif\/ication of $4$-dimensional superbialgebras will be given in a~forthcoming paper.

\subsection{Connected superbialgebras}

In the following we consider connected superbialgebras whose even part are isomorphic to $\mathbb K$.
The notions of a~connected superalgebra, supercoalgebra and Hopf superalgebra are def\/ined in the same way.
\begin{lem}
\label{connexe}
Let $A=A_{0}\oplus A_{1}$ be an $(n+1)$-dimensional connected superspace.
If $(A, \mu, \eta, \Delta, \varepsilon)$ is a~superbialgebra with $1=\eta(1)$, then we have
\begin{enumerate}\itemsep=0pt
\item[$1)$] $\mu(x \otimes y)=\mu(y \otimes x)=0$, $\forall\, x, y \in A_{1}$,
\item[$2)$] $\Delta(x)=1\otimes x + x\otimes 1$, $\forall\, x\in A_{1}$.
\end{enumerate}
\end{lem}

\begin{proof}
Let $(x_{1}, x_{2}, \dots, x_{n})$ be a~basis of $A_{1}$, and set $\mu(x_{i}\otimes x_{j})=\alpha_{ij}1$,
$\Delta(x_{i})=\sum\limits_{k=1}^{n}(\lambda_{k1} x_{k}\otimes 1 +\lambda_{1k} 1\otimes x_{k})$, with
$\alpha_{ij}, \lambda_{1j},\lambda_{j1} \in \mathbb{K}$, $\forall\, i, j \in \{1, \dots, n\}$.
\begin{itemize}\itemsep=0pt
\item The compatibility of the counit and the multiplication
$\varepsilon(\mu(x_{i}\otimes x_{j}))=\varepsilon (x_{i}) \varepsilon(x_{j})$
and the fact that $\varepsilon(x_{i})=0$, $\forall\, i\in \{1, \dots,n\}$,
show that $\alpha_{ij}=0$, $\forall\, i, j \in \{1, \dots, n\}$.
Then we have $\mu(x_{i}\otimes x_{j})=\mu(x_{j}\otimes x_{i})=0$, $\forall\, i, j \in \{1, \dots, n\}$.
\item The compatibility of the counit and the following condition:
\begin{gather*}
(\varepsilon\otimes \id)(\Delta(x_{i}))=(\id\otimes\varepsilon)(\Delta(x_{i}))=x_{i}
,\quad\forall\, i\in\{1,\dots,n\},
\end{gather*}
lead to $(\varepsilon \otimes \id)(\Delta
(x_{i}))=\sum\limits_{k=1}^{n}(\lambda_{k1}\varepsilon(x_{k})\otimes 1+\lambda_{1k} \varepsilon(1)\otimes
x_{k})=\sum\limits_{k=1}^{n}\lambda_{1k} x_{k}=x_{i}$
and
\begin{gather*}
(\id \otimes \varepsilon)(\Delta(x_{i}))
=\sum\limits_{k=1}^{n}(\lambda_{k1} x_{k}\otimes \varepsilon(1)+\lambda_{1k}1\otimes\varepsilon(x_{k}))
=\sum\limits_{k=1}^{n}\lambda_{k1} x_{k}=x_{i}.
\end{gather*}
Then $\lambda_{1k}=\lambda_{k1}=0$, $\forall\, k\in \{1, \dots, n\}\setminus \{i\}$ and $\lambda_{1i}=\lambda_{i1}=1$.
So
\begin{gather*}
\Delta(x_{i})=x_{i}\otimes 1 +1\otimes x_{i},\qquad \forall\, i \in \{1, \dots, n\}.\tag*{\qed}
\end{gather*}
\end{itemize}  \renewcommand{\qed}{}
\end{proof}

Notice that the f\/irst assertion of this lemma is just Proposition~2.11 in~\cite{Armour} for which the
proof is dif\/ferent.
\begin{thm}
\label{$dimA0=1$}
An $n$-dimensional connected superbialgebra exists only when $n < 3$.
\end{thm}
\begin{proof}
Let $A=A_{0}\oplus A_{1}$ be a~superalgebra with multiplication $\mu$ and unit $1$.
Assume $A_{0}=\mathbb K 1$, $\dim A_{1}=n-1 $ and $(x_{1}, x_{2}, \dots, x_{n-1})$ be a~basis of $A_{1}$.

According to Lemma~\ref{connexe}, $\Delta(x_i)=1\otimes x_i + x_i\otimes 1$. Compatibility condition leads
in one hand to $\Delta(\mu(x_{i}\otimes x_{j}))=0$, and in the other hand to
\begin{gather*}
\mu\otimes\mu\circ\tau(\Delta(x_{i})\otimes\Delta(x_{j}))=\mu\otimes\mu\circ\tau((1\otimes x_{i}+x_{i}
\otimes1)\otimes(1\otimes x_{j}+x_{j}\otimes1))
\\
\qquad
=\mu\otimes\mu\circ\tau(x_{i}\otimes1\otimes x_{j}\otimes1+x_{i}\otimes1\otimes1\otimes x_{j}
+1\otimes x_{i}\otimes x_{j}\otimes1+1\otimes x_{i}\otimes1\otimes x_{j})
\\
\qquad
=x_{i}\otimes x_{j}-x_{j}\otimes x_{i},
\qquad
\forall\, i,j\in\{1,\dots,n-1\}.
\end{gather*}
Since $x_{i}\otimes x_{j}-x_{j}\otimes x_{i} \neq 0$, $\forall\, i\neq j $, therefore, the compatibility
condition is satisf\/ied only if $n=2$.
\end{proof}

A $2$-dimensional superbialgebra is either trivial or connected.
For trivial $2$-dimensional bialgebras, we refer to~\cite{D-Makh}.
In the following, we provide connected $2$-dimensional superbialgebras and Hopf superalgebras
classif\/ication.

Let $(A,\mu,\eta,\Delta,\varepsilon)$ be a~$2$-dimensional connected superbialgebra.
We set $A_0=\mathbb K$, $A_1={\rm span}\{x\}$ and $\eta(1)=1$ its unit element.
\begin{prop}
Every $2$-dimensional connected superbialgebra is isomorphic to $2$-dimensional connected superbialgebra
$\mathbb K[x]/(x^2)$ with $\deg(x)=1$ and defined by $\Delta (x)=1\otimes x+x\otimes 1$.

Moreover,
it carries a~Hopf superalgebra structure with an antipode $S$ defined by $S(1)=1$ and $S(x)=-x$.
\end{prop}

\begin{proof}
The multiplication, the comultiplication and the counit are def\/ined according to Lem\-ma~\ref{connexe} and with this extra product, coproduct and counit $x\cdot x=\alpha1$,
$\Delta (x)=\beta 1\otimes x+\gamma x\otimes 1$, where $\alpha,\beta,\gamma \in \mathbb K$.
Solving the system corresponding to conditions~\eqref{asso10}, \eqref{unit10}, with respect to
structure constants $\alpha$, $\beta$, $\gamma$, leads to the result.
For the second assertion, we assume that the antipode $S$ is def\/ined as, $S(1)=1$, $S(x)=\lambda x$, with
$\lambda\in \mathbb K$. Applying the identity~\eqref{antip} to~$x$, we obtain only one non trivial
$2$-dimensional Hopf superalgebra associated to connected $2$-dimensional superbialgebra def\/ined above.
The antipode is def\/ined as $S(1)=1$, $S(x)=-x$.
\end{proof}

Let $f:A\longrightarrow A$ be a~superalgebra homomorphism, this means that $f$ is an even linear map,
satisfying $f(1)=1$ and $ f\circ\mu=\mu\circ f \otimes f$. Let us set $f(x)=\alpha x$, where $\alpha \in
\mathbb K\backslash \{0\}$. A direct calculation shows that $f\circ\mu=\mu\circ f \otimes f$ is
satisf\/ied for any $\alpha \neq 0$.
Therefore, the automorphism group of the $2$-dimensional superbialgebra is the inf\/inite group
\begin{gather*}
\left\{\left(
\begin{matrix}
1&0
\\
0&\alpha
\\
\end{matrix}
\right),
\
\text{with}
\
\alpha\in\mathbb K\backslash\{0\}\right\}.
\end{gather*}

\section[Classification of 3-dimensional superbialgebras and Hopf superalgebras]
{Classif\/ication of 3-dimensional superbialgebras\\ and Hopf superalgebras}\label{classification3}

In dimension $3$, there are three cases for $n_{0}=\dim A_{0}$.
If $n_{0}=1$, we have connected superbialgebras which are covered by Proposition~\ref{$dimA0=1$}.
If $n_{0}=3$, we have trivial superbialgebras for which we refer to~\cite{D-Makh} for the classif\/ication
of $3$-dimensional bialgebras.
It remains to study the case $n_{0}=2$.

\subsection{Superalgebras}

Let $(A, \mu, \eta, \Delta, \varepsilon)$ be a~$3$-dimensional superbialgebra
such that $A=A_{0}\oplus A_{1}$ and $\dim A_{0}=2$.
Let $\{1, x, y\}$ be a~basis of $A$, such that $\{1, x\}$ generates the even part $A_{0}$ and $\{y\}$
generates the odd part $A_{1}$.
Assume that $\eta(1)=1$.

We recall here the classif\/ication of $2$-dimensional algebras.
\begin{prop}[\cite{D-Makh, Gabriel}]\label{pr}
There are, up to isomorphism, two $2$-dimensional algebras $A_1=\mathbb K[x]/(x^2)$
and $A_2=\mathbb K[x]/(x^2-x)$.
\end{prop}

Since the even part of a~superalgebra is an algebra, we f\/ix the even part multiplication to be one of the
two $2$-dimensional algebras, recalled in Proposition~\ref{pr}.
We denote by $\mu_{i|j}$ the $j^{\rm th}$ multiplication obtained by extending the multiplication $\mu_{i}$ of
a~$2$-dimensional algebra to the multiplication of a~$3$-dimensional superalgebra.
We consider f\/irst the algebra $A_{1}$.
\begin{prop}
\label{PropClassDim3-1}
Every non-trivial $3$-dimensional superalgebra, where the even part is the $2$-dimensional algebra $A_{1}$
is isomorphic to one of the following pairwise nonisomorphic $3$-di\-men\-si\-onal superalgebras
$A_{1|1}=\mathbb{K}[x,y]/(x^{2}, y^{2}, xy)$ and $ A_{1|2}=\mathbb{K}[x,y]/(x^{2}, y^{2}-x, xy)$,
where $\deg(x)=0$ and $\deg(y)=1$.
\end{prop}

\begin{proof}
We set
\begin{gather*}
x\cdot y\phantom{}=\alpha y,
\qquad
y\cdot x=\beta y,
\qquad
 y\cdot y=\gamma e_{1}^{0}+\sigma x,
\qquad
\text{with}
\quad
\alpha,\beta,\gamma,\sigma\in\mathbb K.
\end{gather*}
We have $ x \cdot  (y \cdot y )=\gamma x$, $(x \cdot y ) \cdot y=\alpha\gamma
e_{1}^{0}+ \alpha\sigma x$, $(y \cdot y) \cdot x=\gamma x$ and $y \cdot (y \cdot x)=\beta\gamma e_{1}^{0}+ \beta\sigma x$.

According to associativity and by identif\/ication, we obtain $\alpha\gamma=0$, $\alpha\sigma
=\gamma$, $\beta\gamma=0$, $\beta\sigma=\gamma$. Then $\gamma=0$ and $y\cdot y=\sigma x$.

Moreover, $x\cdot (x \cdot y )= (x \cdot x ) \cdot y=0$ and $ (y\cdot x  )
\cdot x=y \cdot  (x \cdot x )=0$ lead to $\alpha=\beta=0$.
Finally, it remains to f\/ix $\sigma$.

If $\sigma=0$, we obtain the superalgebra def\/ined by $A_{1|1}$.
If $\sigma \neq 0$, then we have the superalgebra which is given by the multiplication $\mu_{\sigma}$
def\/ined as $\mu_{\sigma}(x \otimes x)=0$, $\mu_{\sigma}(y \otimes x)=0$, $\mu_{\sigma}(x \otimes
y)=0$, $\mu_{\sigma}(y \otimes y)=\sigma x, \sigma \neq 0.
$ It is isomorphic to $A_{1|2}$ ($\sigma=1$).
\end{proof}

For the second algebra $A_{2}$, we obtain three possible extensions as a~$3$-dimensional superalgebra.
\begin{prop}
\label{PropClassDim3-2}
Every non-trivial $3$-dimensional superalgebra where the even part is the $2$-dimensional algebra $A_{2}$,
is isomorphic to one of the following pairwise nonisomorphic $3$-dimensional superalgebras
$A_{2|1}=\mathbb{K}[x,y]/(x^{2}-x, y^{2}-x,xy-y)$, $A_{2|2}=\mathbb{K}\langle x,y\rangle/(x^{2}-x, y^{2},$ $xy-y, yx)$
and $A_{2|3}=\mathbb{K}[x,y]/(x^{2}-x$, $y^{2},xy-y)$, where $\deg(x)=0$ and $\deg(y)=1$.
\end{prop}

\begin{proof}
By associativity conditions $\left(x \cdot x\right)\cdot y=x \cdot \left(x\cdot y\right)$ and $ y\cdot
\left(x \cdot x\right)=\left(y \cdot x\right)\cdot x$, it follows that $\alpha$ and $\beta$ may have only
values $0$ or $1$.

Similarly, conditions $ x \cdot \left(y \cdot y\right)=\left(x \cdot y \right)\cdot y$ and $ y \cdot
\left(y \cdot x\right)=\left(y \cdot y \right)\cdot x$ lead to
\begin{gather}
\alpha\gamma=0,
\qquad
 \alpha\sigma=\gamma+\sigma,
\qquad
 \beta\gamma=0,
\qquad
 \beta\sigma=\gamma+\sigma.
\label{beta}
\end{gather}
We discuss the two cases $\gamma=0$ and $\gamma\neq 0$.

1.~If $\gamma=0$, then the system~\eqref{beta} reduces to $\sigma\left(\alpha-1\right)=0$
and $\sigma\left(\beta-1\right)=0$.
So, we deduce two subcases:
\begin{enumerate}
\itemsep=0pt
\item[a)] If $\sigma\neq 0$, then $\alpha=\beta=1$.
We obtain superalgebras which are isomorphic to $A_{2|1}$ ($\sigma=1$).

\item[b)] If $\sigma=0$, we study two cases $\alpha=\beta$ and $\alpha\neq\beta$,
\begin{enumerate}\itemsep=0pt
\item[i)] If $\alpha=\beta$, that means $(\alpha,\beta)\in\left\{(0,0),(1,1)\right\}$.
Then, we obtain two superalgebras isomorphic to $A_{2|3}$.

\item[ii)] If $\alpha\neq\beta$, that means $(\alpha, \beta)\in\left\{(0,1),(1,0)\right\}$.
Then, we obtain superalgebras which are isomorphic to $A_{2|2}$.
\end{enumerate}
\end{enumerate}

2.~If $\gamma\neq 0$.
Then, the system~\eqref{beta} reduces to $\alpha=\beta=0$ and $\sigma=-\gamma$.
So, we have superalgebras with multiplications $\mu_{\gamma}$ def\/ined as
\begin{gather*}
\mu_{\gamma}(x\otimes x)=x,
\qquad
 \mu_{\gamma}(y\otimes x)=0,
\qquad
 \mu_{\gamma}(x\otimes y)=0,
\qquad
 \mu_{\gamma}
(y\otimes y)=\gamma\left(1-x\right),
\end{gather*}
which are isomorphic to the superalgebra $A_{2|1}$.
\end{proof}

\subsection{Superbialgebras and Hopf superalgebras}

 Now, we construct superbialgebra structures associated
to the f\/ive superalgebras found above.
\begin{prop}
There is no $3$-dimensional superbialgebra, with $\dim A_{0}=2$, associated to superalgebras $A_{1|1}$ and
$A_{1|2}$.
\end{prop}
\begin{proof}
Assume $\Delta (x)=\alpha 1\otimes 1+\beta 1\otimes x+ \gamma x\otimes 1+\sigma x\otimes x + \delta
y\otimes y$ where $\alpha, \beta, \gamma, \sigma, \delta \in \mathbb K$.
Since we have $x^2=0$, the compatibility condition $\varepsilon(\mu(x\otimes
x))=\varepsilon(x)\varepsilon(x)$, implies $\varepsilon(x)=0$.

In one hand, the condition $(\varepsilon \otimes \id)(\Delta (x))=(\id \otimes \varepsilon)(\Delta(x))=x$
implies that $\alpha=0$ and $\beta=\gamma=1$.
In the other hand, the compatibility condition $\Delta \circ \mu (x\otimes x)=(\mu \otimes \mu)\circ (\id
\otimes \tau \otimes \id)\circ (\Delta \otimes \Delta)(x\otimes x)=0, $ leads to $\beta \gamma=0$.
Therefore we have a~contradiction.
\end{proof}

In the sequel, we consider the superalgebra structures def\/ined by $\mu_{2|1}$, $\mu_{2|2}$ and
$\mu_{2|3}$.
We denote by $A_{i|j}^{k}$ the superbialgebra $(A,\mu_{i|j},\eta,\Delta_{i|j}^{k},\varepsilon_{i|j}^{k})$.
We have $A^{\rm cop}=(A,\mu,\eta,\Delta',\varepsilon)$ where $\Delta'=\Delta\circ \tau$ and $\tau$ is the
superf\/lip.
By direct calculation, we obtain the following pairwise non-isomorphic superbialgebras.
For all of them we have $\Delta_{i|j}^{k}(1)=1\otimes 1$.

For superalgebra $A_{2|1}$, we have
\begin{enumerate}\itemsep=0pt
\item[1)] $A_{2|1}^{1}$ with $\Delta_{2|1}^{1}(x)=1\otimes x+x\otimes 1-x\otimes
x$, $\Delta_{2|1}^{1}(y)=y\otimes 1+ 1\otimes y- y\otimes x$, $\varepsilon_{2|1}^{1}(x)=0$;

\item[2)] $(A_{2|1}^{1})^{\rm cop}$.
\end{enumerate}

For superalgebra $A_{2|2}$, we have
\begin{enumerate}\itemsep=0pt
\item[1)] $A_{2|2}^{1}$ with $\Delta_{2|2}^{1}(x)=1\otimes x+x\otimes 1-x\otimes
x$, $\Delta_{2|2}^{1}(y)=1\otimes y+y\otimes 1-y\otimes x-x\otimes y$, $\varepsilon_{2|2}^{1}(x)=0$;

\item[2)] $A_{2|2}^{2}$ with $\Delta_{2|2}^{2}(x)=1\otimes x+x\otimes 1-x\otimes x+y\otimes
y$, $\Delta_{2|2}^{2}(y)=1\otimes y+y\otimes 1-y\otimes x-x\otimes y$, $\varepsilon_{2|2}^{2}(x)=0$;

\item[3)] $A_{2|2}^{3}$ with $\Delta_{2|2}^{3}(x)=x\otimes x$, $\Delta_{2|2}^{3}(y)=y\otimes x+ x\otimes y$, $\varepsilon_{2|2}^{3}(x)=1$;

\item[4)] $A_{2|2}^{4}$ with $\Delta_{2|2}^{4}(x)=x\otimes x+y\otimes y$, $\Delta_{2|2}^{4}(y)=y\otimes x+
x\otimes y$, $\varepsilon_{2|2}^{4}(x)=1$.
\end{enumerate}

 For superalgebra $A_{2|3}$, we have
\begin{enumerate}\itemsep=0pt
\item[1)] $A_{2|3}^{1}$ with $\Delta_{2|3}^{1}(x)=1\otimes x+ x\otimes 1- x\otimes
x$, $\Delta_{2|3}^{1}(y)=1\otimes y+ y\otimes 1- x\otimes y- y\otimes x$, $\varepsilon_{2|3}^{1}(x)=0$;

\item[2)] $A_{2|3}^{2}$ with $\Delta_{2|3}^{2}(x)=1\otimes x+ x\otimes 1- x\otimes
x$, $\Delta_{2|3}^{2}(y)=1\otimes y+ y\otimes 1- x\otimes y$, $\varepsilon_{2|3}^{2}(x)=0$;

\item[3)] $(A_{2|3}^{2})^{\rm cop}$;

\item[4)] $A_{2|3}^{4}$ with $\Delta_{2|3}^{4}(x)=1\otimes x+ x\otimes 1- x\otimes
x$, $\Delta_{2|3}^{4}(y)=y\otimes 1+ 1\otimes y, \varepsilon_{2|3}^{4}(x)=0$;

\item[5)] $A_{2|3}^{5}$ with $\Delta_{2|3}^{5}(x)=x\otimes x$, $\Delta_{2|3}^{5}(y)=y\otimes x+ x\otimes
y$, $\varepsilon_{2|3}^{5}(x)=1$.
\end{enumerate}

\begin{thm}
Every non-trivial $3$-dimensional superbialgebra with $\dim A_{0}=2$ is isomorphic to one of the following
$3$-dimensional pairwise non-isomorphic superbialgebras, which are defined above,
\begin{gather*}
A_{2|1}^{1},
\quad
\big(A_{2|1}^{1}\big)^{\rm cop},
\quad
A_{2|2}^{k}
\ \ (k=1,\dots,4),
\quad
A_{2|3}^{1},
\quad
 A_{2|3}^{2},
\quad
\big(A_{2|3}^{2}\big)^{\rm cop},
\quad
 A_{2|3}^{4},
\quad
 A_{2|3}^{5}.
\end{gather*}
\end{thm}
\begin{thm}
There is no non-trivial $3$-dimensional Hopf superalgebra.
\end{thm}
\begin{proof}
We check that no one of the superbialgebras may carry a~structure of Hopf superalgebra.
Let $S$ be an antipode of one of the $3$-dimensional superbialgebras def\/ined above.
Assume that
\begin{gather*}
S(1)=1,
\qquad
 S(x)=\lambda_{1}1+\lambda_{2}x,
\qquad
 S(y)=\lambda_{3}y,
\end{gather*}
and $S$ satisf\/ies the identity $\mu \circ (S\otimes \id) \circ \Delta=\mu \circ(\id \otimes S)\circ \Delta
=\eta \circ \varepsilon.
$

We apply the identity to $x$ and study two cases.

{\bf Case 1:} $\varepsilon(x)=0$.
For all $3$-dimensional superalgebras in this case, we have
\begin{gather*}
\Delta(x)=1\otimes x+x\otimes1-x\otimes x+\alpha y\otimes y,
\qquad
\text{such that}
\quad
\alpha=0
\quad
\text{or}
\quad
1.
\end{gather*}
In one hand, we have,
\begin{gather*}
\mu\circ(S\otimes \id_{A})\circ\Delta(x)=\mu\circ(1\otimes x+x\otimes1-S(x)\otimes x+\alpha S(y)\otimes y)
\\
\qquad
=x+(\lambda_{1}x+\lambda_{2}x)-\lambda_{1}x-\lambda_{2}x
=\lambda_{1}1+(1-\lambda_{1})x.
\end{gather*}
In the other hand, we have $\eta \circ \varepsilon(x)=0$. Therefore, we have a~contradiction.

{\bf Case 2:} $\varepsilon(x)=1$.
In this case, $\Delta(x)$ must have the form
\begin{gather*}
\Delta(x)=x\otimes x+\alpha y\otimes y,
\qquad
\text{such that}
\quad
\alpha=0
\quad
\text{or}
\quad
1.
\end{gather*}
Then, on one hand, we have,
\begin{gather*}
\mu\circ(S\otimes \id_{A})\circ\Delta(x)=\mu\circ(S(x)\otimes x+\alpha S(y)\otimes y)=(\lambda_{1}
+\lambda_{2})x.
\end{gather*}
On the other hand, we have $\eta \circ \varepsilon(x)=1$.
So we have a~contradiction.
\end{proof}

\section{Classif\/ication of 4-dimensional superbialgebras\\ and Hopf superalgebras}\label{classification4}

Let $A=A_0\oplus A_1$ be a~$4$-dimensional superbialgebra.
Then one has to consider four cases for $n_{0}=\dim A_{0}$.
If $n_{0}=1$, the superbialgebras are connected, then it is covered by Proposition~\ref{$dimA0=1$}.
If $n_{0}=4$, then they correspond to trivial superbialgebras.
In the sequel, we discuss cases $n_{0}=2$ and $n_{0}=3$.

\subsection[$4$-dimensional algebras]{$\boldsymbol{4}$-dimensional algebras} In the following, we recall
the classif\/ication of $4$-dimensional algebras given by Gabriel in~\cite{Gabriel} and the
classif\/ication of $4$-dimensional superalgebras with $n_{0}=2$ and $n_{0}=3$ provided by Armour, Chen and
Zhang in~\cite{Armour1}.
\begin{thm}
\label{algebra}
The following algebras are pairwise non-isomorphic and every $4$-dimensional algebras is isomorphic to one
of these $19$ algebras:

$1)$ $\mathbb K \times \mathbb K \times \mathbb K \times \mathbb K $,

$2)$ $\mathbb K \times \mathbb K
\times \mathbb K [x]/x^{2}$,

$3)$ $\mathbb K [x]/x^{2} \times \mathbb K [y]/y^{2}$,

$4)$ $\mathbb K \times\mathbb K [x]/x^{3}$,

$5)$ $\mathbb K [x]/x^{4}$,

$6)$ $\mathbb K \times\mathbb K [x,y]/(x,y)^{2}$,

$7)$ $\mathbb K[x,y]/(x^{2}, y^{2})$,

$8)$ $\mathbb K[x,y]/(x^{3}, xy, y^{2})$,

$9)$ $\mathbb K [x, y, z]/(x, y, z)^{2}$,

$10)$ $ M_{2}=\left(
\begin{matrix}
 \mathbb K & \mathbb K
\\
\mathbb K & \mathbb K
\\
\end{matrix}
\right)$,

$11)$ $\left\{ \left(
\begin{matrix}
 a & 0 & 0 & 0
\\
0 & a & 0 & d
\\
c & 0 & b & 0
\\
0 & 0 & 0 & b
\\
\end{matrix}
\right)/ a, b, c, d \in \mathbb K \right\}$,

$12)$ $\Lambda \mathbb K^{2}= \text{exterior algebra of}~\mathbb
K^{2}$,

$13)$ $\mathbb K \times \left(
\begin{matrix}
 \mathbb K & \mathbb K
\\
0 & \mathbb K
\\
\end{matrix}
\right)$,

$14)$ $\left\{ \left(
\begin{matrix}
 a & 0 & 0
\\
c & a & 0
\\
d & 0 & b
\\
\end{matrix}
\right)/ a, b, c, d \in \mathbb K \right\}$,

$15)$ $\left\{ \left(
\begin{matrix}
 a & c & d
\\
0 & a & 0
\\
0 & 0 & b
\\
\end{matrix}
\right)/ a, b, c, d \in \mathbb K \right\}$,

$16)$ $\mathbb K[x,y]/(x^{2}, y^{2}, yx)$,

$17)$ $\left\{ \left(
\begin{matrix}
 a & 0 & 0
\\
0 & a & 0
\\
c & d & b
\\
\end{matrix}
\right)/ a, b, c, d \in \mathbb K \right\}$,

$18)$ $\mathbb K \langle x, y \rangle /(x^{2}, y^{2}, yx-\lambda xy)$, for $\lambda\neq -1, 0, 1$,

$19)$
$\mathbb K \langle x, y \rangle /(y^{2}, x^{2}+yx, xy+yx)$,

\noindent
where $\mathbb K \langle x, y\rangle$ is the free associative algebra generated by $x$ and $y$.
\end{thm}

\subsection[Superalgebras with $\dim(A_{0})=3$]{Superalgebras with $\boldsymbol{\dim(A_{0})=3}$}

 Let
$\{e_{1}^{0}, e_{2}^{0}, e_{3}^{0}, e_{1}^{1}\}$ be a~basis of the superalgebra $A$, $A=A_{0}\oplus
A_{1}$, such that $A_0={\rm span}\{e_{1}^{0}, e_{2}^{0}, e_{3}^{0}\}$, $A_1={\rm span}\{e_{1}^{1}\}$ and $e_{1}^{0}$
be the unit element of the superalgebra.
We recall the multiplication of all possible graduation obtained from a~f\/ixed algebra, we denote by
$\mu_{i|j}$ the multiplication of each $4$-dimensional superalgebra in the case $\dim A_{0}=3$.
We preserve the notation used in~\cite{Armour1}, a~superalgebra $i|j$ means the $j^{\rm th}$ superalgebra
obtained by a~f\/ixed $i^{\rm th}$ $4$-dimensional algebra.

\begin{prop}[\cite{Armour1}]\label{armourprop$1$}
Let $\mathbb K$ be an algebraically closed field of characteristic $0$, suppose that $A$
is a~$4$-dimensional superalgebra with $\dim A_{0}=3$.
Then $A$ is isomorphic to one of the superalgebra in the following pairwise nonisomorphic families:
\begin{description}\itemsep=0pt
\item[${\bf 1|1}$:] $\mathbb K \times \mathbb K \times \mathbb K \times \mathbb K$, $e_{1}^{0}=(1, 1,
1,1)$, $e_{2}^{0}=(1, 0, 0,0)$, $e_{3}^{0}=(0, 0, 1,1)$, $e_{1}^{1}=(0, 0, 1, -1)$;
\item[${\bf 2|1}$:] $\mathbb K \times \mathbb K \times \mathbb K [x]/x^{2}$, $e_{1}^{0}=(1, 1,1)$,
$e_{2}^{0}=(1,0,0)$, $e_{3}^{0}=(0,1,0)$, $e_{1}^{1}=(0, 0, x)$;
\item[${\bf 2|2}$:]$\mathbb K \times
\mathbb K \times \mathbb K [x]/x^{2}$, $e_{1}^{0}=(1, 1,1)$, $e_{2}^{0}=(1, 1,0)$, $e_{3}^{0}=(0, 0,
x)$, $e_{1}^{1}=(1, -1,0)$;
\item[${\bf 3|1}$:] $\mathbb K [x]/x^{2} \times \mathbb K [y]/y^{2}$, $e_{1}^{0}=(1,1)$, $e_{2}^{0}=(1,0)$,
$e_{3}^{0}=(x,0)$, $e_{1}^{1}=(0, y)$;
\item[${\bf 4|1}$:] $\mathbb K
\times\mathbb K [x]/x^{3}$, $e_{1}^{0}=(1,1)$, $e_{2}^{0}=(1,0)$, $e_{3}^{0}=(0, x^{2})$, $e_{1}^{1}=(0,
x)$;
\item[${\bf 6|1}$:] $\mathbb K \times\mathbb K [x,y]/(x,y)^{2}:$ $ e_{1}^{0}=(1,1)$, $e_{2}^{0}=(1,0)$,
$e_{3}^{0}=(0, x)$, $e_{1}^{1}=(0, y)$;
\item[${\bf 7|1}$:] $\mathbb K [x,y]/(x^{2}, y^{2})$, $e_{1}^{0}=1$, $e_{2}^{0}=x+y$, $e_{3}^{0}=xy$,
$e_{1}^{1}=x-y$;
\item[${\bf 8|1}$:] $\mathbb K
[x,y]/(x^{3}, xy, y^{2})$, $e_{1}^{0}=1$, $e_{2}^{0}=x$, $e_{3}^{0}=x^{2}$, $e_{1}^{1}=y$;
\item[${\bf 8|2}$:] $\mathbb K [x,y]/(x^{3}, xy, y^{2})$, $e_{1}^{0}=1$, $e_{2}^{0}=x^{2}$, $e_{3}^{0}=y$,
$e_{1}^{1}=x$;
\item[${\bf 9|1}$:] $\mathbb K [x,y,z]/(x,y,z)^{2}$, $e_{1}^{0}=1$, $e_{2}^{0}=x$, $e_{3}^{0}=y$, $e_{1}^{1}=z$;
\item[${\bf 11|1}$:] $\left\{ {\left(
\begin{matrix}
 a & 0 & 0 & 0
\\
0 & a& 0 & d
\\
c & 0 & b & 0
\\
0 & 0 & 0 & b
\\
\end{matrix}
\right)}/ a, b, c, d \in \mathbb K \right\}$, $e_{1}^{0}={\left(
\begin{matrix}
 1 & 0 & 0 & 0
\\
0 & 1 & 0 & 0
\\
0 & 0 & 1 & 0
\\
0 & 0 & 0 & 1
\\
\end{matrix}
\right)}$, $e_{2}^{0}={\left(
\begin{matrix}
 1 & 0 & 0 & 0
\\
0 & 1 & 0 & 0
\\
0 & 0 & 0 & 0
\\
0 & 0 & 0 & 0
\\
\end{matrix}
\right)},$
\\
$ e_{3}^{0}={\left(
\begin{matrix}
 0 & 0 & 0 & 0
\\
0 & 0 & 0 & 1
\\
0 & 0 & 0 & 0
\\
0 & 0 & 0 & 0
\\
\end{matrix}
\right)}$, $e_{1}^{1}={\left(
\begin{matrix}
 0 & 0 & 0 & 0
\\
0 & 0 & 0 & 0
\\
1 & 0 & 0 & 0
\\
0 & 0 & 0 & 0
\\
\end{matrix}
\right)}$;
\item[${\bf 13|1}$:] $\mathbb K \times \left(
\begin{matrix}
 \mathbb K & \mathbb K
\\
0 & \mathbb K
\\
\end{matrix}
\right)=\left\{\left(a,{\left(
\begin{matrix}
 b & c
\\
0 & d
\end{matrix}
\right)}\right)/ a, b, c, d \in \mathbb K \right\}$, $e_{1}^{0}=\left(1,{\left(
\begin{matrix}
 1 & 0
\\
0 & 1
\\
\end{matrix}
\right)}\right)$,
\\
$ e_{2}^{0}=\left(0,{\left(
\begin{matrix}
 1 & 0
\\
0 & 0
\end{matrix}
\right)}\right)$, $e_{3}^{0}=\left(0,{\left(
\begin{matrix}
 0 & 0
\\
0 & 1
\end{matrix}
\right)}\right)$, $e_{1}^{1}=\left(0,{ \left(
\begin{matrix}
 0 & 1
\\
0 & 0
\end{matrix}
\right)}\right)$;
\item[${\bf 14|1}$:] $\left\{ {\left(
\begin{matrix}
 a & 0 & 0
\\
c & a & 0
\\
d & 0 & b
\end{matrix}
\right)}/ a, b, c, d \in \mathbb K \right\}$, $e_{1}^{0}={\left(
\begin{matrix}
 1 & 0 & 0
\\
0 & 1 & 0
\\
0 & 0 & 1
\end{matrix}
\right)}$, $e_{2}^{0}={\left(
\begin{matrix}
 1 & 0 & 0
\\
0 & 1 & 0
\\
0 & 0 & 0
\end{matrix}
\right)},$
\\
$ e_{3}^{0}={\left(
\begin{matrix}
 0 & 0 & 0
\\
0 & 0 & 0
\\
1 & 0 & 0
\end{matrix}
\right)}$, $e_{1}^{1}={\left(
\begin{matrix}
 0 & 0 & 0
\\
1 & 0 & 0
\\
0 & 0 & 0
\end{matrix}
\right)}$;
\item[${\bf 14|2}$:] $\left\{ {\left(
\begin{matrix}
 a & 0 & 0
\\
c & a & 0
\\
d & 0 & b
\\
\end{matrix}
\right)}/ a, b, c, d \in \mathbb K \right\}$, $e_{1}^{0}={\left(
\begin{matrix}
 1 & 0 & 0
\\
0 & 1 & 0
\\
0 & 0 & 1
\end{matrix}
\right)}$, $e_{2}^{0}={\left(
\begin{matrix}
 1 & 0 & 0
\\
0 & 1 & 0
\\
0 & 0 & 0
\end{matrix}
\right)}$,
\\
$e_{3}^{0}={\left(
\begin{matrix}
 0 & 0 & 0
\\
1 & 0 & 0
\\
0 & 0 & 0
\end{matrix}
\right)}$, $e_{1}^{1}={\left(
\begin{matrix}
 0 & 0 & 0
\\
0 & 0 & 0
\\
1 & 0 & 0
\end{matrix}
\right)}$;
\item[${\bf 15|1}$:] $\left\{ {\left(
\begin{matrix}
 a & c & d
\\
0 & a & 0
\\
0 & 0 & b
\end{matrix}
\right)}/ a, b, c, d \in \mathbb K \right\}$, $e_{1}^{0}={{\left(
\begin{matrix}
 1 & 0 & 0
\\
0 & 1 & 0
\\
0 & 0 & 1
\end{matrix}
\right)}}$, $e_{2}^{0}={\left(
\begin{matrix}
 1 & 0 & 0
\\
0 & 1 & 0
\\
0 & 0 & 0
\end{matrix}
\right)}$,
\\
$ e_{3}^{0}={\left(
\begin{matrix}
 0 & 0 & 1
\\
0 & 0 & 0
\\
0 & 0 & 0
\end{matrix}
\right)}$, $e_{1}^{1}={{\left(
\begin{matrix}
 0 & 1 & 0
\\
0 & 0 & 0
\\
0 & 0 & 0
\end{matrix}
\right)}}$;
\item[${\bf 15|2}$:]$ \left\{{{\left(
\begin{matrix}
 a & c & d
\\
0 & a & 0
\\
0 & 0 & b
\end{matrix}
\right)}}/ a, b, c, d \in \mathbb K \right\}$, $e_{1}^{0}={{\left(
\begin{matrix}
 1 & 0 & 0
\\
0 & 1 & 0
\\
0 & 0 & 1
\end{matrix}
\right)}}$, $e_{2}^{0}={\left(
\begin{matrix}
 1 & 0 & 0
\\
0 & 1 & 0
\\
0 & 0 & 0
\end{matrix}
\right)} $,\\
$ e_{3}^{0}={\left(
\begin{matrix}
 0 & 1 & 0
\\
0 & 0 & 0
\\
0 & 0 & 0
\end{matrix}
\right)}$, $e_{1}^{1}={\left(
\begin{matrix}
 0 & 0 & 1
\\
0 & 0 & 0
\\
0 & 0 & 0
\end{matrix}
\right)}$;
\item[${\bf 17|1}$:]$ \left\{{\left(
\begin{matrix}
 a & 0 & 0
\\
0 & a & 0
\\
c & d & b
\end{matrix}
\right)}/ a, b, c, d \in \mathbb K \right\}$, $e_{1}^{0}={\left(
\begin{matrix}
 1 & 0 & 0
\\
0 & 1 & 0
\\
0 & 0 & 1
\end{matrix}
\right)}$, $e_{2}^{0}={\left(
\begin{matrix}
 1 & 0 & 0
\\
0 & 1 & 0
\\
0 & 0 & 0
\end{matrix}
\right)}$,
\\
$ e_{3}^{0}={\left(
\begin{matrix}
 0 & 0 & 0
\\
0 & 0 & 0
\\
1 & 0 & 0
\end{matrix}
\right)}$, $e_{1}^{1}={\left(
\begin{matrix}
 0 & 0 & 0
\\
0 & 0 & 0
\\
0 & 1 & 0
\end{matrix}
\right)}$.
\end{description}
\end{prop}

\subsection [Superbialgebras and Hopf superalgebras with $\dim(A_{0})=3$]{Superbialgebras and Hopf
superalgebras with $\boldsymbol{\dim(A_{0})=3}$} We look for all possible superbialgebras which could be
obtained from a~f\/ixed superalgebra.
The computation are done using a~computer algebra system.
For the convenience of the presentation, we only summarize the results in the following proposition.
The notations and details are collected in Appendix~\ref{appendix B}.
\begin{prop}
Let $(A$, $\mu$, $\eta$, $\Delta$, $\varepsilon)$ be a~$4$-dimensional superbialgebra with $\dim A_{0}=3$.
Then $A$ is isomorphic to one of the superbialgebra in the following pairwise nonisomorphic families:
\begin{center}
\begin{tabular}{cl}
\hline
{superalgebra}&{associated superbialgebras}
\\
\hline
{${\bf 1|1}$}&{$\big(A,\mu_{1|1},\eta,\Delta_{1|1}^{k},\varepsilon_{1|1}^{k}\big)$, $k=1,\dots,12$\tsep{3pt}\bsep{3pt}}
\\
{${\bf 2|1}$}&{$\big(A,\mu_{2|1},\eta,\Delta_{2|1}^{k},\varepsilon_{2|1}^{k}\big)$, $k=1,\dots,22$\bsep{3pt}}
\\
{${\bf 4|1}$}&{$\big(A,\mu_{4|1},\eta,\Delta_{4|1}^{k},\varepsilon_{4|1}^{k}\big)$, $k=1,2,3$\bsep{3pt}}
\\
{${\bf 6|1}$}&{$\big(A,\mu_{6|1},\eta,\Delta_{6|1}^{k},\varepsilon_{6|1}^{k}\big)$, $k=1,\dots,18$\bsep{3pt}}
\\
{${\bf 13|1}$}&{$\big(A,\mu_{13|1},\eta,\Delta_{13|1}^{k},\varepsilon_{13|1}^{k}\big)$, $k=1,\dots,21$\bsep{3pt}}
\\
{${\bf 14|1}$}&{$\big(A,\mu_{14|1},\eta,\Delta_{14|1}^{k},\varepsilon_{14|1}^{k}\big)$, $k=1,\dots,9$\bsep{3pt}}
\\
{${\bf 14|2}$}&{$\big(A,\mu_{14|2},\eta,\Delta_{14|2}^{k},\varepsilon_{14|2}^{k}\big)$, $k=1,\dots,4$\bsep{3pt}}
\\
{${\bf 15|1}$}&{$\big(A,\mu_{15|1},\eta,\Delta_{15|1}^{k},\varepsilon_{15|1}^{k}\big)$, $k=1,\dots,9$\bsep{3pt}}
\\
{${\bf 15|2}$}&{$\big(A,\mu_{15|2},\eta,\Delta_{15|2}^{k},\varepsilon_{15|2}^{k}\big)$, $k=1,\dots,4$\bsep{3pt}}
\\
{${\bf 17|1}$}&{$\big(A,\mu_{17|1},\eta,\Delta_{17|1}^{k},\varepsilon_{17|1}^{k}\big)$, $k=1,\dots,11$\bsep{3pt}}
\\
\hline
\end{tabular}
\end{center}
\end{prop}

\begin{rem}
There is no $4$-dimensional superbialgebra with $\dim A_{0}=3$ and underlying multiplications $\mu_{2|2}$,
$\mu_{3|1}$, $\mu_{7|1}$, $\mu_{8|1}$, $\mu_{8|2}$, $\mu_{9|1}$, $\mu_{11|1}$.
\end{rem}

Now, we look for Hopf superalgebra structures.
For a~f\/ixed superbialgebra def\/ined above, we add the antipode's property.
It turns out that there exists only one non-trivial $4$-dimensional Hopf superalgebra in the case $\dim
A_{0}=3$.
It corresponds to the algebra $(1|1)$ and comultiplication $\Delta_{1|1}^{2}$, where $A=\mathbb K \times
\mathbb K \times \mathbb K \times \mathbb K$.
We set $ x=e_1^0- 2e_2^0-e_3^0=(-1, 1, 0,0)$ and $y=e_1^1=(0, 0, 1,-1)$ with $\deg(x)=0$,
$\deg(y)=1$.
It leads to basis $\{1,x,x^2,y\}$ and the algebra may be written as $\mathbb K[x,y]/(x^2+y^2-1,xy)$ with
$\deg(x)=0$ and $\deg(y)=1$.

\begin{prop}
\label{Hopf22dim4}
Every non-trivial $4$-dimensional Hopf superalgebra where $\dim A_{0}=3$ is isomorphic to the
$4$-dimensional Hopf superalgebra $\mathbb K[x,y]/(x^2+y^2-1,xy)$ with $\deg(x)=0$, $\deg(y)=1$ and such that
\begin{gather*}
\Delta(x)=x\otimes x-\alpha y\otimes y,
\qquad
 \varepsilon(x)=1,
\qquad
 S(x)=x,
\\
\Delta(y)=x\otimes y+y\otimes x,
\qquad
 \varepsilon(y)=0,
\qquad
 S(y)=\alpha y.
\end{gather*}
where $\alpha$ is a~primitive $4^{th}$ root of unity.

\end{prop}

\subsection [Superalgebras with $\dim(A_{0})=2$]{Superalgebras with $\boldsymbol{\dim(A_{0})=2}$}

Let
$\{e_{1}^{0}$, $e_{2}^{0}$, $e_{1}^{1}$, $e_{2}^{1}\}$ be a~basis of the underlying superspace $A$, such that
$\{e_{1}^{0}$, $e_{2}^{0}\}$ is a~basis of the even part, $\{e_{1}^{1}$, $e_{2}^{1}\}$ a~basis of the odd
part and $e_{1}^{0}$ is the unit of the superalgebra.
\begin{prop}[\cite{Armour1}] \label{armourprop$2$}Let $\mathbb K$ be an algebraically closed f\/ield, suppose that $A$ is a~$4$-dimensional
superalgebra with $\dim A_{0}=2$.
Then $A$ is isomorphic to one of the superalgebra in the following pairwise non-isomorphic families:
\begin{description}\itemsep=0pt
\item[${\bf 1|2}$:] $\mathbb K \times \mathbb K \times \mathbb K \times \mathbb K$,
$e_{1}^{0}=(1, 1, 1,1)$, $e_{2}^{0}=(1, 1, 0,0)$, $e_{1}^{1}=(1, -1, 0,0)$, $e_{2}^{1}=(0, 0, 1, -1)$;
\item[${\bf 2|3}$:] $\mathbb K \times \mathbb K \times \mathbb K [x]/x^{2}$, $e_{1}^{0}=(1, 1,1)$,
$e_{2}^{0}=(1, 1,0)$, $e_{1}^{1}=(1, -1,0)$, $e_{2}^{1}=(0, 0, x)$;
\item[${\bf 3|2}$:] $\mathbb K[x]/x^{2} \times \mathbb K [y]/y^{2}$, $e_{1}^{0}=(1,1)$, $e_{2}^{0}=(1,0)$,
$e_{1}^{1}=(x,0)$, $e_{2}^{1}=(0, y)$;
\item[${\bf 3|3}$:] $\mathbb K[x]/x^{2} \times \mathbb K [y]/y^{2}$, $e_{1}^{0}=(1,1)$, $e_{2}^{0}=(x, y)$,
$e_{1}^{1}=(1, -1)$, $e_{2}^{1}=(x, -y)$;
\item[${\bf 5|1}$:] $\mathbb K [x]/x^{4}$, $e_{1}^{0}=1$, $e_{2}^{0}=x^{2}$, $e_{1}^{1}=x$, $e_{2}^{1}=x^{3}$;
\item[${\bf 6|2}$:] $\mathbb K \times\mathbb K[x,y]/(x,y)^{2}$, $e_{1}^{0}=(1,1)$, $e_{2}^{0}=(1,0)$,
$e_{1}^{1}=(0,x)$, $e_{2}^{1}=(0,y)$;
\item[${\bf 7|2}$:] $\mathbb K[x,y]/(x^{2},y^{2})$, $e_{1}^{0}=1$, $e_{2}^{0}=x$, $e_{1}^{1}=y$, $e_{2}^{1}=xy$;
\item[${\bf 7|3}$:] $\mathbb K[x,y]/(x^{2},y^{2})$, $e_{1}^{0}=1$, $e_{2}^{0}=xy$, $e_{1}^{1}=x$, $e_{2}^{1}=y$;
\item[${\bf 8|3}$:] $\mathbb K[x,y]/(x^{3}, xy, y^{2})$, $e_{1}^{0}=1$, $e_{2}^{0}=x^{2}$, $e_{1}^{1}=x$, $e_{2}^{1}=y$;
\item[${\bf 9|2}$:] $\mathbb K [x, y, z]/(x, y, z)^{2}$, $e_{1}^{0}=1$, $e_{2}^{0}=x$, $e_{1}^{1}=y$, $e_{2}^{1}=z$;
\item[${\bf 10|1}$:] $ M_{2}=\left(
\begin{matrix}
 \mathbb K & \mathbb K
\\
\mathbb K & \mathbb K
\end{matrix}
\right)=\left\{ \left(
\begin{matrix}
 a & b
\\
c & d
\end{matrix}
\right)/ a, b, c, d \in \mathbb K \right\}$, $e_{1}^{0}=\left(
\begin{matrix}
 1 & 0
\\
0 & 1
\end{matrix}
\right)$, $e_{2}^{0}=\left(
\begin{matrix}
 1 & 0
\\
0 & 0
\end{matrix}
\right)$,
\\
$ e_{1}^{1}=\left(
\begin{matrix}
 0 & 1
\\
0 & 0
\end{matrix}
\right)$, $e_{2}^{1}=\left(
\begin{matrix}
 0 & 0
\\
1 & 0
\end{matrix}
\right)$;
\item[${\bf 11|2}$:] $\left\{ \left(
\begin{matrix}
 a & 0 & 0 & 0
\\
0 & a & 0 & d
\\
c & 0 & b & 0
\\
0 & 0 & 0 & b
\end{matrix}
\right) / a, b, c, d \in \mathbb K \right\}$, $e_{1}^{0}= \left(
\begin{matrix}
 1 & 0 & 0 & 0
\\
0 & 1 & 0 & 0
\\
0 & 0 & 1 & 0
\\
0 & 0 & 0 & 1
\end{matrix}
\right)$,
\\
$e_{2}^{0}= \left(
\begin{matrix}
 1 & 0 & 0 & 0
\\
0 & 1 & 0 & 0
\\
0 & 0 & 0 & 0
\\
0 & 0 & 0 & 0
\end{matrix}
\right) $, $e_{1}^{1}= \left(
\begin{matrix}
 0 & 0 & 0 & 0
\\
0 & 0 & 0 & 1
\\
0 & 0 & 0 & 0
\\
0 & 0 & 0 & 0
\end{matrix}
\right) $, $e_{2}^{1}= \left(
\begin{matrix}
 0 & 0 & 0 & 0
\\
0 & 0 & 0 & 0
\\
1 & 0 & 0 & 0
\\
0 & 0 & 0 & 0
\end{matrix}
\right) $;
\item[${\bf 11|3}$:] $\left\{  \left(
\begin{matrix}
 a & 0 & 0 & 0
\\
0 & a & 0 & d
\\
c & 0 & b & 0
\\
0 & 0 & 0 & b
\end{matrix}
\right) / a, b, c, d \in \mathbb K \right\}$, $e_{1}^{0}= \left(
\begin{matrix}
 1 & 0 & 0 & 0
\\
0 & 1 & 0 & 0
\\
0 & 0 & 1 & 0
\\
0 & 0 & 0 & 1
\end{matrix}
\right)$,
\\
$e_{2}^{0}= \left(
\begin{matrix}
 0 & 0 & 0 & 0
\\
0 & 0 & 0 & 1
\\
1 & 0 & 0 & 0
\\
0 & 0 & 0 & 0
\end{matrix}
\right) $, $e_{1}^{1}= \left(
\begin{matrix}
 1 & 0 & 0 & 0
\\
0 & 1 & 0 & 0
\\
0 & 0 & -1 & 0
\\
0 & 0 & 0 & -1
\end{matrix}
\right) $, $e_{2}^{1}= \left(
\begin{matrix}
 0 & 0 & 0 & 0
\\
0 & 0 & 0 & -1
\\
1 & 0 & 0 & 0
\\
0 & 0 & 0 & 0
\end{matrix}
\right) $;
\item[${\bf 12|1}$:] $\Lambda \mathbb K^{2}=\mathbb K \langle x, y \rangle / (x^{2}, y^{2},
xy+yx)$, $e_{1}^{0}=1$, $e_{2}^{0}=x$, $e_{1}^{1}=y$, $e_{2}^{1}=xy$;

\item[${\bf 12|2}$:] $\Lambda
\mathbb K^{2}=\mathbb K \langle x, y \rangle / (x^{2}, y^{2}, xy+yx)$, $e_{1}^{0}=1$, $e_{2}^{0}=xy$, $e_{1}^{1}=x$, $e_{2}^{1}=y$;
\item[${\bf 14|3}$:] $ \left\{ \left(
\begin{matrix}
 a & 0 & 0
\\
c & a & 0
\\
d & 0 & b
\end{matrix}
\right) / a, b, c, d \in \mathbb K \right\}$, $e_{1}^{0}= \left(
\begin{matrix}
 1 & 0 & 0
\\
0 & 1 & 0
\\
0 & 0 & 1
\end{matrix}
\right) $, $e_{2}^{0}= \left(
\begin{matrix}
 1 & 0 & 0
\\
0 & 1 & 0
\\
0 & 0 & 0
\end{matrix}
\right) $,
\\
$e_{1}^{1}= \left(
\begin{matrix}
 0 & 0 & 0
\\
1 & 0 & 0
\\
0 & 0 & 0
\end{matrix}
\right) $, $e_{2}^{1}= \left(
\begin{matrix}
 0 & 0 & 0
\\
0 & 0 & 0
\\
1 & 0 & 0
\end{matrix}
\right) $;
\item[${\bf 15|3}$:] $ \left\{  \left(
\begin{matrix}
 a & c & d
\\
0 & a & 0
\\
0 & 0 & b
\end{matrix}
\right) / a, b, c, d \in \mathbb K \right\}$, $e_{1}^{0}= \left(
\begin{matrix}
 1 & 0 & 0
\\
0 & 1 & 0
\\
0 & 0 & 1
\end{matrix}
\right) $, $e_{2}^{0}= \left(
\begin{matrix}
 1 & 0 & 0
\\
0 & 1 & 0
\\
0 & 0 & 0
\end{matrix}
\right) $,
\\
$e_{1}^{1}= \left(
\begin{matrix}
 0 & 1 & 0
\\
0 & 0 & 0
\\
0 & 0 & 0
\end{matrix}
\right) $, $e_{2}^{1}= \left(
\begin{matrix}
 0 & 0 & 1
\\
0 & 0 & 0
\\
0 & 0 & 0
\end{matrix}
\right) $;

\item[${\bf 16|1}$:] $\mathbb K \langle x, y \rangle / (x^{2}, y^{2}, yx)$,
$e_{1}^{0}=1$, $e_{2}^{0}=x$, $e_{1}^{1}=y$, $e_{2}^{1}=x y$;

\item[${\bf 16|2}$:] $\mathbb K \langle x, y
\rangle / (x^{2},y^{2},yx), $ $e_{1}^{0}=1$, $e_{2}^{0}=y$, $e_{1}^{1}=x$, $e_{2}^{1}=x y$;

\item[${\bf 16|3}$:] $\mathbb K \langle x, y\rangle / (x^{2}, y^{2}, yx)$, $e_{1}^{0}=1$, $e_{2}^{0}=x
y$, $e_{1}^{1}=x$, $e_{2}^{1}=y$;

\item[${\bf 17|2}$:] $\left\{ \left(
\begin{matrix}
 a & 0 & 0
\\
0 & a & 0
\\
c & d & b
\end{matrix}
\right) / a, b, c, d \in \mathbb K \right\}$, $e_{1}^{0}= \left(
\begin{matrix}
 1& 0 & 0
\\
0 & 1 & 0
\\
0 & 0 & 1
\end{matrix}
\right) $, $e_{2}^{0}= \left(
\begin{matrix}
 1 & 0 & 0
\\
0 & 1 & 0
\\
0 & 0 & 0
\end{matrix}
\right)$,
\\
$e_{1}^{1}= \left(
\begin{matrix}
 0 & 0 & 0
\\
0 & 0 & 0
\\
1 & 0 & 0
\end{matrix}
\right)$, $e_{2}^{1}= \left(
\begin{matrix}
 0 & 0 & 0
\\
0 & 0 & 0
\\
0 & 1 & 0
\end{matrix}
\right)$;
\item[${\bf (18;\lambda|1)}$:] $\mathbb K \langle x, y \rangle / (x^{2}, y^{2}, yx-\lambda
xy)$, where $\lambda \in \mathbb K $ with $\lambda \neq -1, 0, 1 $,
\\
$ e_{1}^{0}=1$, $e_{2}^{0}=x$, $e_{1}^{1}=y$, $e_{2}^{1}=x y$;

\item[${\bf (18;\lambda|2)}$:] $\mathbb K \langle x, y \rangle / (x^{2}, y^{2}, yx-\lambda xy)$, where $\lambda \in \mathbb K$ with $\lambda \neq -1, 0, 1 $,
\\
$ e_{1}^{0}=1$, $e_{2}^{0}=xy$, $e_{1}^{1}=x$, $e_{2}^{1}=y$;

\item[${\bf 19|1}$:]$\mathbb K \langle x,
y\rangle /(y^{2}, x^{2}+yx, xy+yx)$, $e_{1}^{0}=1$, $e_{2}^{0}=x y$, $e_{1}^{1}=x$, $e_{2}^{1}=y$.
\end{description}
\end{prop}

\subsection[Superbialgebras and Hopf superalgebras with $\dim(A_{0})=2$]{Superbialgebras and Hopf
superalgebras with $\boldsymbol{\dim(A_{0})=2}$}

Now, we compute all superbialgebra structures associated
to superalgebras described above where $\dim A_{0}=2$.
We give fairly basic description of the results, all details are given in Appendix~\ref{appendix C}.
We denote by $\mu_{i|j}$ the multiplication in $A_{i|j}$ with respect to the basis described above.
\begin{prop}
There is no $4$-dimensional superbialgebra with $\dim A_{0}=2$ endowed with one of these multiplications
$\mu_{1|2}$, $\mu_{3|3}$, $\mu_{5|1}$, $\mu_{7|2}$, $\mu_{7|3}$, $\mu_{8|3}$, $\mu_{10|1}$, $\mu_{11|3}$, $\mu_{12|1}$,
$\mu_{16|1}$,
$ \mu_{16|2}$, $\mu_{16|3}$, $\mu_{18|1}$, $\mu_{18|2}$, $\mu_{19|1}$.
\end{prop}

\begin{prop}
Let $(A$, $\mu$, $\eta$, $\Delta$, $\varepsilon)$ be a~$4$-dimensional superbialgebra with $\dim A_{0}=2$, then
$A$ is isomorphic to one of the superbialgebra in the following pairwise nonisomorphic families:

\begin{center}
\begin{tabular}{cl}
\hline
{superalgebra}&{associated superbialgebras}
\\
\hline
{${\bf 2|3}$ }&{$\big(A, \mu_{2|3}, \eta, \Delta_{2|3}^{k},\varepsilon_{2|3}^{k}\big)$, $k=1, \dots, 4$\tsep{3pt}\bsep{3pt}}
\\
{${\bf 3|2}$ }&{$\big(A, \mu_{3|2}, \eta, \Delta_{3|2}^{k},\varepsilon_{3|2}^{k}\big)$, $k=1, \dots, 9$\bsep{3pt}}
\\
{${\bf 6|2}$}&{$\big(A, \mu_{6|2}, \eta, \Delta_{6|2}^{k},\varepsilon_{6|2}^{k}\big)$, $k=1, \dots, 11$}
\\
{${\bf 11|2}$}&{$\big(A, \mu_{11|2}, \eta, \Delta_{11|2}^{1},\varepsilon_{11|2}^{1}\big)$\bsep{3pt}}
\\
{${\bf 12|2}$}&{$\big(A, \mu_{12|2}, \eta, \Delta_{12|2}^{1},\varepsilon_{12|2}^{1}\big)$\bsep{3pt}}
\\
{${\bf 14|3}$}&{$\big(A, \mu_{14|3}, \eta, \Delta_{14|3}^{1},\varepsilon_{14|3}^{1}\big)$, $k=1, \dots, 7$\bsep{3pt}}
\\
{${\bf 15|3}$}&{$\big(A, \mu_{15|3}, \eta, \Delta_{15|3}^{1},\varepsilon_{15|3}^{1}\big)$, $k=1, \dots, 7$\bsep{3pt}}
\\
{${\bf 17|2}$}&{$ \big(A, \mu_{17|2}, \eta, \Delta_{17|2}^{1},\varepsilon_{17|2}^{1}\big)$, $k=1, 2$\bsep{3pt} }
\\
\hline
\end{tabular}
\end{center}
\end{prop}

It turns out that the only superbialgebras, with $\dim A_{0}=2$, which can have a~non-trivial Hopf
superalgebra structure are $(3|2)$, $(11|2)$ and $(12|2)$.
\begin{prop}
\label{alg}
Every non-trivial $4$-dimensional Hopf superalgebra where $\dim A_{0}=2$ is isomorphic to one of the
following pairwise non-isomorphic Hopf superalgebras:
\begin{enumerate}\itemsep=0pt
\item[$1)$] $\mathcal{H}_{1}=\mathbb K\langle x,y\rangle /( x^{2}-x,y^{2},yx)$ with $\deg(x)=0$,
$\deg(y)=1$ and such that
\begin{gather*}
\Delta(x)=1\otimes x+x\otimes1-2x\otimes x,
\qquad
 \Delta(y)=1\otimes y+y\otimes1,
\\
\varepsilon(x)=\varepsilon(y)=0,
\qquad
 S(x)=x,
\qquad
 S(y)=-y;
\end{gather*}
\item[$2)$] $\mathcal{H}_{2}=\mathbb K\langle x,y\rangle /( x^{2}-x,y^{2},xy-yx-y) $ with $\deg(x)=0$,
$\deg(y)=1$ and such that
\begin{gather*}
\Delta(x)=1\otimes x+x\otimes1-2x\otimes x,
\qquad
 \Delta(y)=1\otimes y+y\otimes1-2x\otimes xy-2xy\otimes x,
\\
\varepsilon(x)=\varepsilon(y)=0,
\qquad
 S(x)=x,
\qquad
 S(y)=y;
\end{gather*}

\item[$3)$] $\mathcal{H}_{3}=\mathbb K\langle x,y\rangle /( x^{2},y^{2},xy+yx) $ with $\deg(x)=\deg(y)=1$ and
such that
\begin{gather*}
\Delta(x)=1\otimes x+x\otimes1,
\qquad
 \Delta(y)=1\otimes y+y\otimes1,
\\
\varepsilon(x)=\varepsilon(y)=0,
\qquad
 S(x)=-x,
\qquad
 S(y)=-y;
\end{gather*}
\item[$4)$] $\mathcal{H}_{4}=\mathbb K\langle x,y\rangle /( x^{2}-x,y^{2},yx)$ with $\deg(x)=0$, $\deg(y)=1$ and
such that
\begin{gather*}
\Delta(x)=1\otimes x+x\otimes1-2x\otimes x,
\qquad
 \Delta(y)=1\otimes y+y\otimes1-2x\otimes y,
\\
\varepsilon(x)=\varepsilon(y)=0,
\qquad
 S(x)=x,
\qquad
 S(y)=2xy-y,
\end{gather*}
\end{enumerate}
where $\mathbb K\langle x,y\rangle $ stands for noncommutative polynomials.
\end{prop}

\begin{proof}
The Hopf superalgebra $\mathcal{H}_{1}$ corresponds to $A_{3|2}^{2}$ (see Appendix~\ref{appendix C}) with
$x=e_{2}^{0}$, $y=e_{1}^{1}+e_{2}^{1}$.

The Hopf superalgebra $\mathcal{H}_{2}$ corresponds to $A_{11|2}^{1}$ with $x=e_{2}^{0}$, $y=e_{1}^{1}-e_{2}^{1}$.

The Hopf superalgebra $\mathcal{H}_{3}$ corresponds to $A_{12|2}^{1}$ with $x=e_{1}^{1}$, $y=e_{2}^{1}$.

The Hopf superalgebra $\mathcal{H}_{4}$ corresponds to $A_{3|2}^{1}$ with $x=e_{2}^{0}$, $y=e_{1}^{1}+e_{2}^{1}$.
\end{proof}

Therefore, gathering the results corresponding to the case $dim A_0=3$ and $dim A_0=2$ we obtain f\/ive
$4$-dimensional Hopf superalgebras.

\begin{thm}
Every non-trivial $4$-dimensional Hopf superalgebra is isomorphic to one of the following Hopf superalgebras
\begin{gather*}
\mathcal{H}_1\cong\mathbb K[\mathbb{Z}/2\mathbb{Z}]\otimes\Lambda\mathbb{K},
\qquad
\mathcal{H}_2\cong\mathbb K[\mathbb{Z}/2\mathbb{Z}]\rtimes_\sigma\Lambda\mathbb{K},
\\
\mathcal{H}_3\cong\Lambda\mathbb{K}^2,
\qquad
\mathcal{H}_4=\mathcal{H}_2^*
\qquad
\text{and}
\qquad
\mathcal{H}_5,
\end{gather*}
where $\sigma:\mathbb{Z}/2\mathbb{Z}\rightarrow GL(\mathbb K )$ is the non-trivial action of
$\mathbb{Z}/2\mathbb{Z}$ on $\mathbb K$, $\mathcal{H}_2^*$ is the dual of $\mathcal{H}_2$ and
$\mathcal{H}_5$ is defined as $\mathbb K[x,y]/(x^2+y^2-1,xy)$ $(\deg(x)=0$, $\deg(y)=1)$ such that
\begin{gather*}
\Delta(x)=x\otimes x-\alpha y\otimes y,
\qquad
 \varepsilon(x)=1,
\qquad
 S(x)=x,
\\
\Delta(y)=x\otimes y+y\otimes x,
\qquad
 \varepsilon(y)=0,
\qquad
 S(y)=\alpha y.
\end{gather*}
with $\alpha^4=1$.
\end{thm}

\begin{proof}
The f\/irst three Hopf superalgebras are cocommutative.
We obtain the isomorphisms, thanks to Kostant's theorem.
However, since $\mathcal{H}_4$ is commutative, its dual is cocommutative, it turns out that we have
$\mathcal{H}_4=\mathcal{H}_2^*$.
The Hopf superalgebra $\mathcal{H}_5$ corresponds to $A_{1|1}$ (Proposition~\ref{Hopf22dim4}), it is not
cocommutative and is the only one whose even part is $3$-dimensional.
\end{proof}

In the following table, we collect the results obtained for $4$-dimensional superbialgebras and Hopf
superalgebras.
Notice that algebras are those classif\/ied by Gabriel, see Theorem~\ref{algebra}.
The superalgebras $ (i|j)$ denote $j^{\rm th}$ graduation of the $i^{\rm th}$ algebra, obtained by Armour, Chen and
Zhang, see Proposition~\ref{armourprop$1$} for the f\/irst case ($\dim A_{0}=3$) and
Proposition~\ref{armourprop$2$} for the second case ($\dim A_{0}=2$).
We quote below the number of corresponding superbialgebras and Hopf superalgebras.

\begin{center}
\begin{tabular}
{c c c c }
\hline
{algebra}&{  superalgebra}&{$\sharp$ superbialgebras}&{$\sharp$ Hopf
superalgebras}
\\
\hline
{$1$}&{${1|1}$ }&{ 12 }&{ 1 }
\\
{}&{${1|2}$}&{ 0 }&{ 0 }
\\
\hline  {${2} $}&{${2|1}$}&{ 22 }&{ 0 }
\\
{} &{${2|2}$}&{ 0 }&{ 0 }
\\
{}&{${2|3}$ }&{ 4 }&{ 0 }
\\
\hline {${3}$}&{${3|1}$}&{ 0 }&{ 0 }
\\
{}&{${3|2}$}&{ 9 }&{ 2 }
\\
{}&{${3|3}$}&{ 0 }&{ 0 }
\\
\hline {${4} $}&{${4|1}$}&{ 3 }&{ 0 }
\\
\hline {${5}$}&{${5|1}$}&{ 0 }&{ 0 }
\\
\hline {${6}$}&{${6|1}$ }&{ 18 }&{ 0 }
\\
{} &{${6|2}$ }&{ 11 }&{ 0 }
\\
\hline {${7} $}&{${7|1}$}&{ 0 }&{ 0 }
\\
{} &{${7|2}$ }&{ 0 }&{ 0 }
\\
{}&{${7|3}$}&{ 0 }&{ 0 }
\\
\hline {${8}$}&{${8|1}$}&{ 0 }&{ 0 }
\\
{} &{${8|2}$}&{ 0 }&{ 0 }
\\
{}&{${8|3}$}&{ 0 }&{ 0 }
\\
\hline {${9}$}&{${9|1}$}&{ 0 }&{ 0 }
\\
{} &{${9|2}$}&{ 0 }&{ 0 }
\\
\hline {${10}$}&{${10|1}$}&{ 0 }&{ 0 }
\\
\hline {${11}$}&{${11|1}$}&{ 0 }&{ 0 }
\\
{} &{${11|2}$}&{ 1 }&{ 1 }
\\
{}&{${11|3}$}&{ 0 }&{ 0 }
\\
\hline {${12}$}&{${12|1}$}&{ 0 }&{ 0 }
\\
{} &{${12|2}$ }&{ 1 }&{ 1 }
\\
\hline {${13}$}&{${13|1}$}&{ 21 }&{ 0 }
\\
\hline
{${14}$}&{${14|1}$ }&{ 9}&{ 0 }
\\
  {}&{${14|2} $}&{ 4 }&{ 0 }
\\
{}&{${14|3}$}&{ 7 }&{ 0 }
\\
\hline {${15}$}&{${15|1}$}&{ 9 }&{ 0 }
\\
{} &{${15|2} $}&{ 4 }&{ 0 }
\\
{}&{${15|3}$}&{ 7 }&{ 0 }
\\
\hline {${16}$}&{${16|2}$}&{ 0 }&{ 0 }
\\
{} &{${16|1}$}&{ 0 }&{ 0 }
\\
{}&{${16|3}$}&{\;0 }&{ 0 }
\\
\hline {${17}$}
&{${17|1}$}&{ 11 }&{ 0 }
\\
{} &{${17|2}$}&{ 2 }&{ 0 }
\\
\hline {$({18}; \lambda)$}&{$({18}; \lambda)|1$}&{ 0 }&{ 0 }
\\
 {} &{$({18};\lambda)|2$}&{ 0 }&{ 0 }
\\
\hline {${19}$}&{${19|1}$}&{ 0 }&{ 0 }
\\
\hline
\end{tabular}
\end{center}

\appendix

\section{Appendix}

We list in the following supercoalgebras associated to a~given $4$-dimensional superalgebra $A$, such that
$A=A_{0}\oplus A_{1}$.
We denote the comultiplication by $\Delta_{i|j}^{k}$ and the counit by $\varepsilon_{i|j}^{k}$, where~$i$
indicates the item of the $4$-dimensional algebra listed in Theorem~\ref{algebra}, and $A_{i|j}$ denotes
the superalgebra obtained from the $i^{\rm th}$ algebra, see Proposition~\ref{armourprop$1$} for the f\/irst
case ($\dim A_{0}=3$) and Proposition~\ref{armourprop$2$} for the second case ($\dim A_{0}=2$).
The exponent $k$ indicates the item of the comultiplication and counit which combined with the
multiplication of superalgebra $A_{i|j}$ and the unit $\eta$ provide a~$4$-dimensional superbialgebra.
Recall that for all of them we have $\eta(1)=1$ (unit element), $\Delta(1)=1\otimes 1$,
$\varepsilon(1)=1$ and $\varepsilon(A_{1})=0$.
In the sequel, we denote by $A_{i|j}^{k}$ the superbialgebra
$(A,\mu_{i|j},\eta,\Delta_{i|j}^{k},\varepsilon_{i|j}^{k})$.
In order to simplify the superbialgebra in this appendix, we change the variables for superalgebras
$A_{1|1}$, $A_{2|1}$, $A_{3|2}$, $A_{6|1}$, $A_{13|1}$, $A_{11|2}$, $A_{12|2}$.
For the superalgebras $A_{1|1}$, $A_{3|2}$, $A_{11|2}$, $A_{12|2}$ we change the variables as mentioned
before Proposition~\ref{Hopf22dim4} and in the proof of Proposition~\ref{alg}.
For $A_{2|1}$ and $A_{13|1}$ we use $x=e_{2}^{0}-e_{3}^{0}$, $y=e_{1}^{1}$, and for $A_{6|1}$ we use
$x=e_{2}^{0}+e_{3}^{0}$, $y=e_{1}^{1}$.
According to the multiplication of each superalgebra we obtain that $A_{2|1}\cong \mathbb{K}[x,y]/(x^{3}-x, xy,y^{2})$,
$A_{6|1}\cong \mathbb{K}[x,y]/(x^{3}-x^{2}, xy,y^{2})$,
$A_{13|1}\cong\mathbb{K}\langle x, y\rangle /(x^{3}-x, xy+yx,xy-y, y^{2})$.
For the remaining superalgebras we change just the notation of basis vectors.
For $\dim(A_{0})=3$, we consider the basis $\{1=e_{1}^{0}, x=e_{2}^{0}, y=e_{3}^{0}, z=e_{1}^{1}\}$ and
for $\dim(A_{0})=2$ the basis $\{1=e_{1}^{0}, x=e_{2}^{0}, y=e_{1}^{1}, z=e_{2}^{1}\}$.

\subsection[Case $\dim A_{0}=3$]{Case $\boldsymbol{\dim A_{0}=3}$}\label{appendix B}

\textit{Superalgebra $A_{1|1}\cong \mathbb{K}[x,y]/(x^{2}+y^{2}-1, xy)$} with $\deg(x)=0$ and $\deg(y)=1 $,
we have
\begin{enumerate}\itemsep=0pt
\item[1)] $A_{1|1}^{1}$ with $\Delta_{1|1}^{1}(x)=\frac{1}{2}(x \otimes 1+x^{2}\otimes
1+x^{2}\otimes x -x\otimes x),
\\
\Delta_{1|1}^{1}(y)=y\otimes 1+\frac{1}{2}x^{2}\otimes y-\frac{1}{2} x\otimes
y$, $\varepsilon_{1|1}^{1}(x)=-1$;
\item[2)] $A_{1|1}^{2}$ with $\Delta_{1|1}^{2}(x)=x\otimes x-\alpha
y\otimes y$, $\Delta_{1|1}^{2}(y)=x\otimes y+y\otimes x$, $\varepsilon_{1|1}^{2}(x)=1,$
\\
where $\alpha$ is a~primitive 4th root of unity of $\mathbb K$;

\item[3)] $A_{1|1}^{3}$ with $\Delta_{1|1}^{3}(x)=x\otimes x$, $\Delta_{1|1}^{3}(y)=1\otimes y+y\otimes
x$, $\varepsilon_{1|1}^{3}(x)=1$;

\item[4)] $A_{1|1}^{4}$ with $\Delta_{1|1}^{4}(x)=\frac{1}{2}( x^{2}\otimes
x^{2}+x\otimes x^{2}+x^{2}\otimes x-x\otimes x )$,
\\
$\Delta_{1|1}^{4}(y)=1\otimes y+y\otimes x^{2}$, $\varepsilon_{1|1}^{4}(x)=-1$;

\item[5)] $A_{1|1}^{5}$ with $\Delta_{1|1}^{5}(x)=\frac{1}{2}(x\otimes 1+ 1\otimes x-1\otimes x^{2}-x^{2}\otimes 1+ x^{2}\otimes x^{2}
+x\otimes x)$,
\\
$\Delta_{1|1}^{5}(y)=y\otimes 1+\frac{1}{2}(y\otimes x+ x\otimes y-y\otimes x^{2}+x^{2}\otimes y)$,
$\varepsilon_{1|1}^{5}(x)=1$;

\item[6)] $A_{1|1}^{6}$ with $\Delta_{1|1}^{6}(x)=\frac{1}{2}(x\otimes x+ x\otimes x^{2}+x^{2}\otimes
x-x^{2}\otimes x^{2})$,
\\
$\Delta_{1|1}^{6}(y)=1\otimes y+y\otimes x^{2}$, $\varepsilon_{1|1}^{6}(x)=1$;

\item[7)] $A_{1|1}^{7}$ with
$\Delta_{1|1}^{7}(x)=-1\otimes 1+\frac{1}{2}( x\otimes x+x\otimes 1+1\otimes x+1\otimes x^{2}+x^{2}\otimes
1-x^{2}\otimes x^{2})$,
\\
$\Delta_{1|1}^{7}(y)=\frac{1}{2}(x\otimes y+y\otimes x+x^{2}\otimes y+y\otimes
x^{2})$, $\varepsilon_{1|1}^{7}(x)=1$;

\item[8)] $A_{1|1}^{8}$ with $\Delta_{1|1}^{8}(x)=x\otimes
x$, $\Delta_{1|1}^{8}(y)=y\otimes 1+x^{2}\otimes y$, $\varepsilon_{1|1}^{8}(x)=1$;

\item[9)]
$A_{1|1}^{9}=(A_{1|1}^{4})^{\rm cop}$, $A_{1|1}^{10}=(A_{1|1}^{1})^{\rm cop}$, $A_{1|1}^{11}=(A_{1|1}^{3})^{\rm cop}$, $A_{1|1}^{12}=(A_{1|1}^{5})^{\rm cop}$.
\end{enumerate}

\textit{Superalgebra} $A_{2|1}\cong \mathbb{K}[x,y]/(x^{3}-x, y^{2}, xy)$ with $\deg(x)=0$ and $\deg(y)
=1 $,
$ \varepsilon_{2|1}^{k}(x)=1$ for $k=1,\dots, 22$:
\begin{enumerate}\itemsep=0pt
\item[1)] $A_{2|1}^{k}$ with  $\Delta_{2|1}^{k}(x)=x \otimes x$ for $ k=1, \dots, 5$:
\begin{itemize}\itemsep=0pt
\item $\Delta_{2|1}^{1}(y)=y\otimes x +x\otimes y$,
\item $\Delta_{2|1}^{2}(y)=x\otimes y +y\otimes x^{2}$,
\item $\Delta_{2|1}^{3}(y)=y\otimes x^{2} +x^{2}\otimes y$,
\item $\Delta_{2|1}^{4}(y)=1\otimes y +y\otimes x$,
\item $\Delta_{2|1}^{5}(y)=1\otimes y +y\otimes x^{2}$;
\end{itemize}
\item[2)] $A_{2|1}^{6}$ with $\Delta_{2|1}^{6}(x)=\frac{1}{2}(x\otimes 1+ x\otimes x+ x^{2}\otimes x -
x^{2}\otimes 1),
\Delta_{2|1}^{6}(y)=\frac{1}{2}(y\otimes x+ x\otimes y+ x^{2}\otimes y + y\otimes x^{2})$;

\item[3)] $A_{2|1}^{7}$ with $\Delta_{2|1}^{7}(x)=\frac{1}{4}(3 x\otimes x + x\otimes x^{2}+ x^{2}\otimes x -
x^{2}\otimes x^{2})$,
\\
$\Delta_{2|1}^{7}(y)=\frac{1}{2}y\otimes x+\frac{1}{2}x\otimes y+\frac{1}{2}x^{2}\otimes y
+\frac{1}{2}y\otimes x^{2}$;

\item[4)] $A_{2|1}^{8}$ with $\Delta_{2|1}^{8}(x)=\frac{1}{4} (3x\otimes x+
x\otimes x^{2}+ x^{2}\otimes x - x^{2}\otimes x^{2})$,
\\
$\Delta_{2|1}^{8}(y)=1\otimes y+\frac{1}{2} y\otimes x+\frac{1}{2} y\otimes x^{2}$;

\item[5)] $A_{2|1}^{9}$ with
$\Delta_{2|1}^{9}(x)=\frac{1}{2} (\frac{3}{2} x^{2}\otimes x^{2} - \frac{1}{2} x^{2}\otimes x- \frac{1}{2}
x\otimes x^{2} + \frac{3}{2} x\otimes x+ 1\otimes x- x\otimes 1- x^{2}\otimes 1- 1\otimes
x^{2})$, $\Delta_{2|1}^{9}(y)=\frac{1}{2}( y\otimes x+ x\otimes y+ x^{2}\otimes y + y\otimes x^{2})$;

\item[6)]
$A_{2|1}^{10}$ with $\Delta_{2|1}^{10}(x)=\frac{1}{2}( x^{2}\otimes x^{2}+ x\otimes x+ 1\otimes x +
x\otimes 1- 1\otimes x^{2}- x^{2}\otimes 1)$,
\\
$Ê \Delta_{2|1}^{10}(y)=\frac{1}{2}(y\otimes x+ x\otimes y+ x^{2}\otimes y + y\otimes x^{2})$;

\item[7)]
$A_{2|1}^{11}$ with $\Delta_{2|1}^{11}(x)=\frac{1}{2}( x^{2}\otimes x^{2}+ x\otimes x+ 1\otimes x +
x\otimes 1- 1\otimes x^{2} - x^{2}\otimes 1)$,
\\
$\Delta_{(2|1)}^{11}(y)=1\otimes y+\frac{1}{2} (y\otimes x+ x\otimes y+ y\otimes x^{2}- x^{2}\otimes y)$;

\item[8)] $A_{2|1}^{12}$ with $\Delta_{2|1}^{12}(x)=\frac{1}{2} (x\otimes x^{2}+ x\otimes x+ 1\otimes x -
1\otimes x^{2})$,
\\
$\Delta_{2|1}^{12}(y)=1\otimes y +\frac{1}{2} (y\otimes x^{2}+ y\otimes x)$;

\item[9)] $A_{2|1}^{13}$ with $\Delta_{2|1}^{13}(x)=-1\otimes 1+\frac{1}{2}( x\otimes x+x\otimes 1+1\otimes x+1\otimes x^{2}+x^{2}\otimes
1-x^{2}\otimes x^{2})$,
\\
$Ê\Delta_{2|1}^{13}(y)=\frac{1}{2}(y\otimes x+ x\otimes y+ x^{2}\otimes y + y\otimes x^{2})$;

\item[10)] $A_{2|1}^{14}$ with $\Delta_{2|1}^{14}(x)=1 \otimes x+ x\otimes 1- x \otimes x - x^{2}\otimes x -
x\otimes x^{2}$,
\\
$\Delta_{2|1}^{14}(y)=y\otimes 1+ 1\otimes y- x^{2}\otimes y - y\otimes x^{2}$;

\item[11)] $A_{2|1}^{15}$ with $\Delta_{2|1}^{15}(x)=1\otimes x+x\otimes 1 - x^{2}\otimes x$, $\Delta_{2|1}^{15}(y)=y\otimes 1+ 1\otimes
y- x^{2}\otimes y - y\otimes x^{2}$;

\item[12)] $A_{2|1}^{16}$ with $\Delta_{2|1}^{16}(x)=1\otimes x+x\otimes
1+\frac{1}{2}( - x\otimes x - x\otimes x^{2} - x^{2}\otimes x + x^{2}\otimes x^{2})$,
\\
$\Delta_{2|1}^{16}(y)=y\otimes 1+ 1\otimes y- x^{2}\otimes y - y\otimes x^{2}$;

\item[13)] $A_{2|1}^{17}$ with $\Delta_{2|1}^{17}(x)
=1\otimes x + x\otimes 1 - x\otimes x^{2} - x^{2}\otimes x + x^{2}\otimes x^{2}$,
\\
$\Delta_{2|1}^{17}(y)=y\otimes 1+ 1\otimes y- x^{2}\otimes y - y\otimes x^{2}$;

\item[14)] $A_{2|1}^{18}$ with $\Delta_{2|1}^{18}(x)
=-1\otimes 1+1\otimes x^{2}+x^{2}\otimes 1+\frac{1}{2}(x\otimes x+x\otimes x^{2}+x^{2}\otimes x- 3 x^{2}\otimes x^{2})$,
\\
$\Delta_{2|1}^{18}(y)=\frac{1}{2}y\otimes x+\frac{1}{2}x\otimes y+\frac{1}{2}x^{2}\otimes y
+\frac{1}{2}y\otimes x^{2}$;

\item[15)] $A_{2|1}^{19}$ with $\Delta_{2|1}^{19}(x)=\frac{1}{2}(x\otimes
x+x\otimes x^{2}+x^{2}\otimes x - x^{2}\otimes x^{2})$, $\Delta_{2|1}^{19}(y)=1\otimes y +y\otimes x^{2}$;

\item[16)] $A_{2|1}^{20}=(A_{2|1}^{15})^{\rm cop}$, $A_{2|1}^{21}=(A_{2|1}^{2})^{\rm cop}$, $A_{2|1}^{22}=(A_{2|1}^{6})^{\rm cop}$.
\end{enumerate}

\textit{Superalgebra} $A_{4|1}$, we have $\Delta_{4|1}^{k}(x)=x \otimes x$, $\varepsilon_{4|1}^{k}(x)=1$,
$\varepsilon_{4|1}^{k}(y)=0$ for $k=1,\dots, 3$:
\begin{enumerate}\itemsep=0pt
\item[1)] $A_{4|1}^{1}$ with $\Delta_{4|1}^{1}(y)=x\otimes y+y\otimes
x$, $\Delta_{4|1}^{1}(z)=x\otimes z+z\otimes x$;

\item[2)] $A_{4|1}^{2}$ with $\Delta_{4|1}^{2}(y)=x\otimes
y+y\otimes 1$, $\Delta_{4|1}^{2}(z)=z\otimes 1+x\otimes z$;
\item[3)] $A_{4|1}^{3}=(A_{4|1}^{2})^{\rm cop}$.
\end{enumerate}

\textit{Superalgebra} $A_{6|1} \cong \mathbb{K}[x,y]/(x^{3}-x^{2}, y^{2}-xy)$ with $\deg(x)=0$ and
$\deg(y)=1 $, we have $\varepsilon_{6|1}^{k}(x)=1$ for $k=1,\dots, 18$:
\begin{enumerate}\itemsep=0pt
\item[1)] $A_{6|1}^{1}$ with $\Delta_{6|1}^{1}(x)=x \otimes x^{2}+x^{2}\otimes x-x^{2}\otimes
x^{2}$, $\Delta_{6|1}^{1}(y)=y\otimes x^{2} +x^{2}\otimes y$;

\item[2)] $A_{6|1}^{2}$ with $\Delta_{6|1}^{2}(x)=x \otimes x^{2}+x^{2}\otimes x-x^{2}\otimes x^{2}+y\otimes y$, $\Delta_{6|1}^{2}(y)=y\otimes x^{2} +x^{2}\otimes y$;

\item[3)] $A_{6|1}^{k}$ with $\Delta_{6|1}^{k}(x)=x
\otimes x$ for $k=3,\dots,8$ and
\begin{itemize}\itemsep=0pt
\item $\Delta_{6|1}^{3}(y)=x\otimes y +y\otimes x$,
\item $\Delta_{6|1}^{4}(y)=y\otimes x^{2} +x\otimes y$,
\item $\Delta_{6|1}^{5}(y)=y\otimes x^{2} +x^{2}\otimes y$,
\item $\Delta_{6|1}^{6}(y)=y\otimes 1+x\otimes y$,
\item $\Delta_{6|1}^{7}(y)=y\otimes 1+x^{2}\otimes y$,
\item $\Delta_{6|1}^{8}(y)=1\otimes y+y\otimes 1$;
\end{itemize}

\item[4)] $A_{6|1}^{9}$ with $\Delta_{6|1}^{9}(x)=x\otimes x+y\otimes y$, $\Delta_{6|1}^{9}(y)=x\otimes y
+y\otimes x$;

\item[5)] $A_{6|1}^{10}$ with $\Delta_{6|1}^{10}(x)=x \otimes 1-x^{2}\otimes 1 +x^{2}\otimes
x$, $\Delta_{6|1}^{10}(y)=y\otimes x^{2} +x^{2}\otimes y$;

\item[6)] $A_{6|1}^{11}$ with
$\Delta_{6|1}^{11}(x)=x \otimes x^{2}+x^{2}\otimes x-x^{2}\otimes x^{2}$, $\Delta_{6|1}^{11}(y)=y\otimes
1+x^{2}\otimes y$;

\item[7)] $A_{6|1}^{12}$ with $\Delta_{6|1}^{12}(x)=x\otimes 1-x^{2}\otimes
1+x^{2}\otimes x$, $\Delta_{6|1}^{12}(y)=x^{2}\otimes y+y\otimes 1$;

\item[8)]
$A_{6|1}^{13}=(A_{6|1}^{4})^{\rm cop}$, $A_{6|1}^{14}=(A_{6|1}^{11})^{\rm cop}$, $A_{6|1}^{15}=(A_{6|1}^{6})^{\rm cop}$, $A_{6|1}^{16}=(A_{6|1}^{5})^{\rm cop}$,
$A_{6|1}^{17}=(A_{6|1}^{12})^{\rm cop}$, $A_{6|1}^{18}=(A_{6|1}^{10})^{\rm cop}$.
\end{enumerate}

\textit{Superalgebra} $A_{13|1}\cong \mathbb{K}\langle x, y\rangle /(x^{3}-x, y^{2}, xy+yx, xy-y)$
with $\deg(x)=0$ and $\deg(y)=1 $, we have $\varepsilon_{13|1}^{k}(x)=0$ for $k=1,\dots,11$ and
$ \varepsilon_{13|1}^{k}(x)=1$ for $k=12,\dots,21$:
\begin{enumerate}\itemsep=0pt
\item[1)] $A_{13|1}^{1}$ with $\Delta_{13|1}^{1}(x)=1 \otimes x+ x\otimes 1- x^{2}\otimes x$, $\Delta_{13|1}^{1}(y)=1\otimes y +y\otimes 1-x^{2}\otimes y-y\otimes x^{2}$;

\item[2)] $A_{13|1}^{2}$ with $\Delta_{13|1}^{2}(x)=1 \otimes x+ x\otimes 1+ \frac{1}{2}x\otimes x- \frac{1}{2}x\otimes x^{2}-
\frac{1}{2} x^{2}\otimes x- \frac{1}{2} x^{2}\otimes x^{2}$,
\\
$\Delta_{13|1}^{2}(y)=1\otimes y +y\otimes 1-x^{2}\otimes y-y\otimes x^{2}$;

\item[3)] $A_{13|1}^{3}$ with $\Delta_{13|1}^{3}(x)=1 \otimes x+ x\otimes 1+ \frac{1}{2}x\otimes x- \frac{1}{2}x\otimes x^{2}-
\frac{1}{2} x^{2}\otimes x- \frac{1}{2} x^{2}\otimes x^{2}$,
\\
$Ê \Delta_{13|1}^{3}(y)=1\otimes y +y\otimes 1$;

\item[4)] $A_{13|1}^{4}$ with $\Delta_{13|1}^{4}(x)=1
\otimes x+ x\otimes 1+ - x^{2}\otimes x- x^{2}\otimes x^{2}$,
\\
$\Delta_{13|1}^{4}y)=1 \otimes y+y\otimes 1-y\otimes x^{2}-x^{2}\otimes y$;

\item[5)] $A_{13|1}^{5}$ with $\Delta_{13|1}^{5}(x)=1 \otimes x+ x\otimes 1- \frac{1}{2}x\otimes x- \frac{1}{2}x\otimes x^{2}-
\frac{1}{2} x^{2}\otimes x+ \frac{1}{2} x^{2}\otimes x^{2}$,
\\
$\Delta_{13|1}^{5}(y)=1 \otimes y+ y\otimes 1- \frac{1}{2}x\otimes y- \frac{1}{2} y\otimes x- \frac{1}{2}
x^{2}\otimes y-\frac{1}{2} y\otimes x^{2}$;

\item[6)] $A_{13|1}^{6}$ with $\Delta_{13|1}^{6}(x)=1 \otimes
x+ x\otimes 1- \frac{1}{2}x\otimes x- \frac{1}{2}x\otimes x^{2}- \frac{1}{2} x^{2}\otimes x+ \frac{1}{2}
x^{2}\otimes x^{2}-2 y\otimes y$, $\Delta_{13|1}^{6}(y)=1 \otimes y+ y\otimes 1- \frac{1}{2}x\otimes y-
\frac{1}{2} y\otimes x- \frac{1}{2} x^{2}\otimes y-\frac{1}{2} y\otimes x^{2}$;

\item[7)] $A_{13|1}^{7}$ with $\Delta_{13|1}^{7}(x)=1 \otimes x+ x\otimes
1-\frac{1}{2}x\otimes x- \frac{1}{2}x\otimes x^{2}+\frac{1}{2} x^{2}\otimes x^{2}- x^{2}\otimes x$,
\\
$Ê \Delta_{13|1}^{7}(y)=1 \otimes y+ y\otimes 1- \frac{1}{2}x\otimes y- \frac{1}{2} x^{2}\otimes y- y\otimes
x^{2}$;

\item[8)] $A_{13|1}^{8}$ with $\Delta_{13|1}^{8}(x)=1 \otimes x+ x\otimes 1- x\otimes x^{2}$,
$\Delta_{13|1}^{8}(y)=1 \otimes y+ y\otimes 1- x^{2}\otimes y- y\otimes x^{2}$;

\item[9)] $A_{13|1}^{9}$ with
$\Delta_{13|1}^{9}(x)=1 \otimes x+ x\otimes 1- x\otimes x^{2}$, $\Delta_{13|1}^{9}(y)=y \otimes
1+1\otimes y-y\otimes x^{2}$;

\item[10)] $A_{13|1}^{10}=(A_{13|1}^{9})^{\rm cop}$;

\item[11)] $A_{13|1}^{11}$ with $\Delta_{13|1}^{11}(x)=1 \otimes x+ x\otimes 1- x\otimes x- x\otimes x^{2}- x^{2}\otimes x$,
\\
$Ê\Delta_{13|1}^{11}(y)=1 \otimes y+ y\otimes 1- x^{2}\otimes y- y\otimes x^{2}$;

\item[12)] $A_{13|1}^{12}$ with
$\Delta_{13|1}^{12}(x)=\frac{1}{2}1 \otimes x+\frac{1}{2}x \otimes 1-\frac{1}{2} 1 \otimes
x^{2}-\frac{1}{2}x^{2} \otimes 1-\frac{1}{4}x \otimes x^{2}-\frac{1}{4}x^{2} \otimes x+\frac{3}{4}x \otimes
x+\frac{3}{4}x^{2} \otimes x^{2}$, $\Delta_{13|1}^{12}(y)=\frac{1}{2}x \otimes y+\frac{1}{2}y \otimes
x+\frac{1}{2}x^{2} \otimes y+\frac{1}{2}y \otimes x^{2}$;

\item[13)] $A_{13|1}^{k}$ with $\Delta_{13|1}^{k}(y)
=\frac{1}{2}x \otimes y+\frac{1}{2}y \otimes x+\frac{1}{2}x^{2} \otimes y+\frac{1}{2}y\otimes x^{2}$ for $k=13,\dots,21$:
\begin{itemize}\itemsep=0pt
\item $\Delta_{13|1}^{13}(x)=x \otimes x$,
\item $\Delta_{13|1}^{14}(x)=\frac{1}{2}x \otimes x+\frac{1}{2} x \otimes x^{2}+\frac{1}{2}x^{2}
\otimes x-\frac{1}{2}x^{2} \otimes x^{2}+2 y \otimes y$,
\item $\Delta_{13|1}^{15}(x)=\frac{1}{2}x
\otimes x+\frac{1}{2}x^{2} \otimes x+\frac{1}{2} x \otimes 1-\frac{1}{2}x^{2} \otimes 1$,
\item $\Delta_{13|1}^{16}(x)=\frac{1}{2}1 \otimes x+\frac{1}{2}x \otimes 1-\frac{1}{2} 1 \otimes
x^{2}-\frac{1}{2}x^{2} \otimes 1+\frac{1}{2}x \otimes x+\frac{1}{2}x^{2} \otimes x^{2}$,
\item $\Delta_{13|1}^{17}(x)=\frac{3}{4}x \otimes x+\frac{1}{4} x \otimes x^{2}+\frac{1}{4}x^{2} \otimes
x-\frac{1}{4}x^{2} \otimes x^{2}$,
\item $\Delta_{13|1}^{18}(x)=-1 \otimes 1+\frac{1}{2}x \otimes x+ 1
\otimes x^{2}+x^{2} \otimes 1+\frac{1}{2}x \otimes x^{2}+\frac{1}{2}x^{2} \otimes x-\frac{3}{2}x^{2}
\otimes x^{2}$,
\item $\Delta_{13|1}^{19}(x)=-1 \otimes 1+\frac{1}{2}x \otimes x+ 1 \otimes x^{2}+x^{2}
\otimes 1+\frac{1}{2}x \otimes x^{2}+\frac{1}{2}x^{2} \otimes x-\frac{3}{2}x^{2} \otimes x^{2}+2y \otimes
y$,
\item $\Delta_{13|1}^{20}(x)=-1 \otimes 1+\frac{1}{2}x \otimes x+ \frac{1}{2}1 \otimes
x^{2}+\frac{1}{2}x^{2} \otimes 1+\frac{1}{2}x \otimes x^{2}+\frac{1}{2}x^{2} \otimes x-\frac{3}{2}x^{2}
\otimes x^{2}$;
\end{itemize}
\item[14)] $A_{13|1}^{21}=(A_{13|1}^{15})^{\rm cop}$.
\end{enumerate}

\textit{Superalgebra} $A_{14|1}$, we have $\varepsilon_{14|1}^{k}(x)=0$, $\varepsilon_{14|1}^{k}(y)=0$
for $k=1,\dots, 9$:
\begin{enumerate}\itemsep=0pt
\item[1)] $A_{14|1}^{k}$ with $\Delta_{14|1}^{1}(x)=1\otimes x+ x\otimes 1-x\otimes x$ for
$k=1,\dots, 6$:
\begin{itemize}\itemsep=0pt
\item $\Delta_{14|1}^{1}(y)=1\otimes y+y\otimes 1-x\otimes y-y\otimes x$, $\Delta_{14|1}^{1}(z)=1\otimes
z+z\otimes 1$,
\item $\Delta_{14|1}^{2}(y)=1\otimes y+y\otimes 1-x\otimes y-y\otimes x$,
$\Delta_{14|1}^{2}(z)=1\otimes z+z\otimes 1- x\otimes z$,
\item $\Delta_{14|1}^{3}(y)=1\otimes y+y\otimes 1-x\otimes y-y\otimes x+ y\otimes y$,
$\Delta_{14|1}^{3}(z)=1\otimes z+z\otimes 1-x\otimes z$,
\item $\Delta_{14|1}^{4}(y)=1\otimes y+y\otimes 1-x\otimes y-y\otimes x$,
$\Delta_{14|1}^{4}(z)=1\otimes z+z\otimes 1-x\otimes z-z\otimes x$,
\item $\Delta_{14|1}^{5}(y)=1\otimes y+y\otimes 1-x\otimes
y-y\otimes x+ y\otimes y$, $\Delta_{14|1}^{5}(z)=1\otimes z+z\otimes 1-x\otimes z-z\otimes x+z\otimes y$,
\item $\Delta_{14|1}^{6}(y)=1\otimes y+y\otimes 1-x\otimes z-z\otimes x+ y\otimes y$,
$\Delta_{14|1}^{6}(z)=1\otimes z+z\otimes 1-x\otimes z-z\otimes x$;
\end{itemize}
\item[2)] $A_{14|1}^{7}$ with $\Delta_{14|1}^{7}(x)=1\otimes x+ x\otimes 1-x\otimes x+ y\otimes y$,
\\
$Ê\Delta_{14|1}^{7}(y)=1\otimes y+y\otimes 1-x\otimes y-y\otimes x+ y\otimes y$,
$\Delta_{14|1}^{7}(z)=1\otimes z+z\otimes 1$;
\item[3)]
$A_{14|1}^{8}=(A_{14|1}^{2})^{\rm cop}$, $A_{14|1}^{9}=(A_{14|1}^{3})^{\rm cop}$.
\end{enumerate}

\textit{Superalgebra} $A_{14|2}$, we have $\varepsilon_{14|2}^{k}(x)=0$
and $\varepsilon_{14|2}^{k}(y)=0$ for $k=1,\dots, 4$:
\begin{enumerate}\itemsep=0pt
\item[1)] $A_{14|2}^{1}$ with $\Delta_{14|2}^{1}(x)=1\otimes x+x\otimes 1-x\otimes x$,
$\Delta_{14|2}^{1}(y)=1\otimes y+y\otimes 1-y\otimes z-z\otimes x$,
\\
$Ê\Delta_{14|2}^{1}(z)=1\otimes z+z\otimes 1-z\otimes x- x\otimes z$;

\item[2)] $A_{14|2}^{2}$ with $\Delta_{14|2}^{2}(x)=1\otimes x+x\otimes 1-x\otimes x$,
$\Delta_{14|2}^{2}(y)=1\otimes y+y\otimes
1-x\otimes y$,
\\
$Ê\Delta_{14|2}^{2}(z)=1\otimes z+z\otimes 1-x\otimes z-z\otimes x$;

\item[3)] $A_{14|2}^{3}$ with
$\Delta_{14|2}^{3}(x)=1\otimes x+x\otimes 1-x\otimes x$, $\Delta_{14|2}^{3}(y)=1\otimes y+y\otimes
1-x\otimes y-y\otimes x+ y\otimes y$, $\Delta_{14|2}^{3}(z)=1\otimes z+z\otimes 1-z\otimes x- x\otimes z$;

\item[4)] $A_{14|2}^{4}=(A_{14|2}^{2})^{\rm cop}$.
\end{enumerate}

\textit{Superalgebra} $A_{15|1}$, we have
$A_{15|1}^{1}=(A_{14|1}^{1})^{\rm op}$,
$A_{15|1}^{2}=(A_{14|1}^{2})^{\rm op}$,
$A_{15|1}^{3}=(A_{14|1}^{5})^{\rm op}$,
$A_{15|1}^{4}=(A_{14|1}^{9})^{\rm op}$,
$A_{15|1}^{5}=(A_{14|1}^{6})^{\rm op}$,
$A_{15|1}^{6}=(A_{14|1}^{3})^{\rm op}$,
$A_{15|1}^{7}=(A_{14|1}^{4})^{\rm op}$,
$A_{15|1}^{8}=(A_{14|1}^{7})^{\rm op}$,
$A_{15|1}^{9}=(A_{14|1}^{4})^{\rm op,cop}$.

\textit{Superalgebra $A_{15|2}$}, we have
$A_{15|2}^{1}=(A_{14|2}^{1})^{\rm op}$, $A_{15|2}^{2}=(A_{14|2}^{4})^{\rm op}$,
$A_{15|2}^{3}=(A_{14|2}^{2})^{\rm op}$, $A_{15|2}^{4}=(A_{14|2}^{3})^{\rm op}$.

\textit{Superalgebra} $A_{17|1}$, we have $\varepsilon_{17|1}^{k}(x)=1$
and $\varepsilon_{17|1}^{k}(y)=0$ for $k=1,\dots,11$:
\begin{enumerate}\itemsep=0pt
\item[1)] $A_{17|1}^{1}$ with $\Delta_{17|1}^{1}(x)=x\otimes x+y\otimes y+z\otimes z$,
$\Delta_{17|1}^{1}(y)=x\otimes y+y\otimes x+z\otimes z$,
\\
$Ê \Delta_{17|1}^{1}(z)=x\otimes z+z\otimes x+z\otimes y+y\otimes z$;

\item[2)] $A_{17|1}^{2}$ with $\Delta_{17|1}^{2}(x)=x\otimes x+z\otimes z$, $\Delta_{17|1}^{2}(y)=x\otimes y+y\otimes x$,
$\Delta_{17|1}^{2}(z)=z\otimes x+x\otimes z$;

\item[3)] $A_{17|1}^{3}$ with $\Delta_{17|1}^{3}(x)=x\otimes
x+z\otimes z$, $\Delta_{17|1}^{3}(y)=x\otimes y+y\otimes x+ z\otimes z-y\otimes y$,
\\
$\Delta_{17|1}^{3}(z)=z\otimes x+x\otimes z$;

\item[4)] $A_{17|1}^{4}$ with $\Delta_{17|1}^{4}(x)=x\otimes
x+y\otimes y$, $\Delta_{17|1}^{4}(y)=x\otimes y+y\otimes x$,
\\
$\Delta_{17|1}^{4}(z)=z\otimes x+x\otimes z-z\otimes y-y\otimes z$;

\item[5)] $A_{17|1}^{k}$, we have
$\Delta_{17|1}^{k}(x)=x\otimes x$ and $\Delta_{17|1}^{k}(y)=x\otimes y+y\otimes x+y\otimes y$ for $k=5, 6, 7$:
\begin{itemize}\itemsep=0pt
\item $\Delta_{17|1}^{5}(z)=x\otimes z+z\otimes x$,
\item $\Delta_{17|1}^{6}(z)=x\otimes z+z\otimes
x+z\otimes y-y\otimes z$,
\item $\Delta_{17|1}^{7}(z)=z\otimes x+x\otimes z+z\otimes y$;
\end{itemize}

\item[6)] $A_{17|1}^{8}=(A_{17|1}^{7})^{\rm cop}$;

\item[7)] $A_{17|1}^{k}$, we have $\Delta_{17|1}^{k}(x)=x\otimes
x$ for $k=9, 10, 11$:
\begin{itemize}\itemsep=0pt
\item $\Delta_{17|1}^{9}(y)=x\otimes y+y\otimes x-y\otimes y$, $\Delta_{17|1}^{9}(z)=z\otimes x+x\otimes z$,
\item $\Delta_{17|1}^{10}(y)=x\otimes y+y\otimes x+z\otimes z$, $\Delta_{17|1}^{10}(z)=x\otimes z+z\otimes x$,
\item $\Delta_{17|1}^{11}(y)=x\otimes y+y\otimes x+y\otimes y+z\otimes z$,
$\Delta_{17|1}^{11}(z)=z\otimes x+x\otimes z+z\otimes y+y\otimes z$.
\end{itemize}
\end{enumerate}

\subsection[Case $\dim A_{0}=2$]{Case $\boldsymbol{\dim A_{0}=2}$}
\label{appendix C}

\textit{Superalgebra} $A_{2|3}$, we have $\Delta_{2|3}^{k}(x)=1\otimes x+x\otimes 1-x\otimes x$,
$\varepsilon_{2|3}^{k}(x)=0$ for $k=1,\dots,4$:
\begin{enumerate}\itemsep=0pt
\item[1)] $A_{2|3}^{1}$ with $\Delta_{2|3}^{1}(y)=1\otimes y+ y\otimes 1-x\otimes
y$, $\Delta_{2|3}^{1}(z)=1\otimes z+z\otimes 1-x\otimes z-z\otimes x$;

\item[2)] $A_{2|3}^{2}$ with $
\Delta_{2|3}^{2}(y)=1\otimes y+ y\otimes 1-x\otimes y+z\otimes x$, $\Delta_{2|3}^{2}(z)=1\otimes z+z\otimes
1-x\otimes z-z\otimes x$;

\item[3)] $A_{2|3}^{3}=(A_{2|3}^{1})^{\rm cop}$, $A_{2|3}^{4}=(A_{2|3}^{2})^{\rm cop}$.
\end{enumerate}

\textit{Superalgebra $A_{3|2}\cong \mathbb{ K}\langle x, y \rangle /(x^{2} - x, y^{2}, yx)$} with $\deg(x)
=0$ and $\deg(y)=1$, we have $\varepsilon_{3|2}^{k}(x)=0$ for $ k=1, \dots, 9$:
\begin{enumerate}\itemsep=0pt
\item[1)] $A_{3|2}^{k}$ with $\Delta_{3|2}^{k}(x)=1\otimes x+x\otimes1-2x\otimes x$ for $k=1, 2$,
and
\begin{itemize}\itemsep=0pt
\item $\Delta_{3|2}^{1}(y)=1\otimes y+y\otimes 1-2x\otimes y$,
\item $\Delta_{3|2}^{2}(y)=1\otimes y+y\otimes 1$;
\end{itemize}

\item[2)] $A_{3|2}^{k}$ with $\Delta_{3|2}^{k}(x)=1\otimes x+x\otimes 1- x\otimes x$, for $k=3, \dots, 9$, and
\begin{itemize}\itemsep=0pt
\item $\Delta_{3|2}^{3}(y)=1\otimes y+y\otimes 1- x\otimes y- y\otimes x+ xy\otimes x+x\otimes xy$,
\item $\Delta_{3|2}^{4}(y)=1\otimes y+y\otimes 1- x\otimes y- xy\otimes x$,
\item $\Delta_{3|2}^{5}(y)=1\otimes y+y\otimes 1- x\otimes y- y\otimes x$,
\item $\Delta_{3|2}^{6}(y)=1\otimes y+y\otimes 1- x\otimes xy-xy\otimes x$,
\item $\Delta_{3|2}^{7}(y)=1\otimes y+y\otimes 1- y\otimes x$,
\item $\Delta_{3|2}^{8}(y)=1\otimes y+y\otimes 1$,
\item $A_{3|2}^{9}=(A_{3|2}^{4})^{\rm cop}$.
\end{itemize}
\end{enumerate}

\textit{Superalgebra} $A_{6|2}$, we have $\Delta_{6|2}^{k}(x)=x\otimes x$ and $\varepsilon_{6|2}^{k}(x)=1$
for $k=1,\dots, 11$:
\begin{enumerate}\itemsep=0pt
\item[1)] $A_{6|2}^{k}$ with $\Delta_{6|2}^{k}(y)=y \otimes x+x\otimes y$ for $k=1,\dots,4$, and
\begin{itemize}\itemsep=0pt
\item $\Delta_{6|2}^{1}(z)=z \otimes 1+1\otimes z$,
\item $\Delta_{6|2}^{2}(z)=z \otimes x+ x\otimes z$,
\item $\Delta_{6|2}^{3}(z)=z \otimes x+1\otimes z$,
\item $\Delta_{6|2}^{4}(z)=1 \otimes z+z\otimes x+1 \otimes y- x\otimes y$;
\end{itemize}
\item[2)] $A_{6|2}^{5}$ with $\Delta_{6|2}^{6}(y)=y \otimes 1+1\otimes y$, $\Delta_{6|2}^{5}(z)=z \otimes x+
x\otimes z$;

\item[3)] $A_{6|2}^{6}$ with $\Delta_{6|2}^{6}(y)=y \otimes x+1\otimes y$, $\Delta_{6|2}^{6}(z)=z
\otimes x+1\otimes z$;

\item[4)] $A_{6|2}^{7}$ with $\Delta_{6|2}^{7}(y)=1 \otimes y+y\otimes
x$, $\Delta_{6|2}^{7}(z)=x \otimes z+z\otimes x$;

\item[5)] $A_{6|2}^{8}$ with $\Delta_{6|2}^{8}(y)=y \otimes
x+x\otimes y+z\otimes 1+1\otimes z-x\otimes z-z\otimes x$, $\Delta_{6|2}^{8}(z)=1\otimes z+z\otimes 1$;

\item[6)] $A_{6|2}^{9}$ with $\Delta_{6|2}^{9}(y)=y \otimes x+1\otimes y+z\otimes 1+1\otimes z-z\otimes
x-x\otimes z$, $\Delta_{6|2}^{9}(z)=z \otimes x+1\otimes z$;

\item[7)] $A_{6|2}^{10}$ with
$\Delta_{6|2}^{10}(y)=y \otimes x+\frac{1}{2} x\otimes y+\frac{1}{2} 1 \otimes y+\frac{1}{4} 1\otimes
z-\frac{1}{4} x \otimes z$,
\\
$\Delta_{6|2}^{10}(z)=z \otimes x+\frac{1}{2} 1\otimes z+1 \otimes y-x\otimes y+\frac{1}{2} x\otimes z$;

\item[8)] $A_{6|2}^{11}=(A_{3|2}^{6})^{\rm cop}$.
\end{enumerate}

\textit{Superalgebra $A_{11|2}\cong \mathbb{K}\langle x,y\rangle /(x^{2} -x,y^{2}, xy-yx-y)$} with $\deg(x)=0, \deg(y)=1$:
\begin{enumerate}\itemsep=0pt
\item[1)] $\Delta_{11|2}^{1}(x)=1\otimes x+x\otimes1-2x\otimes x$, $\Delta_{11|2}^{1}(y)=1\otimes
y+y\otimes 1-2x\otimes xy-2xy\otimes x$,
\\
$\varepsilon_{11|2}^{1}(x)=1$.
\end{enumerate}

\textit{Superalgebra $A_{12|2}^{1}\cong \mathbb{K}\langle x, y\rangle /(x^{2}, y^{2}, xy + yx)$} with
$\deg(x)=\deg(y)=1$:
\begin{enumerate}\itemsep=0pt
\item[1)] $\Delta_{12|2}^{1}(x)=1\otimes x+x\otimes 1$, $\Delta_{12|2}^{1}(y)=1\otimes y+y\otimes
1$.
\end{enumerate}

\textit{Superalgebra $A_{14|3}$}, we have $\Delta_{14|3}^{k}(z)=1 \otimes z+ z\otimes 1- x\otimes z-
z\otimes x$ and $\varepsilon_{14|3}^{k}(x)=0$ for $k=1,\dots,7$:
\begin{enumerate}\itemsep=0pt
\item[1)] $A_{14|3}^{k}$ with $\Delta_{14|3}^{k}(x)=1 \otimes x+ x\otimes 1- x\otimes x$ for
$k=1,\dots,4$:
\begin{itemize}\itemsep=0pt
\item $\Delta_{14|3}^{1}(y)=1 \otimes y+ y\otimes 1$,
\item $\Delta_{14|3}^{2}(y)=1\otimes y+y\otimes1+
x\otimes z+z\otimes x$,
\item $\Delta_{14|3}^{3}(y)=1\otimes y+ y\otimes 1- x\otimes y$,
\item
$\Delta_{14|3}^{4}(y)=1 \otimes y+ y\otimes 1- x\otimes y-y\otimes x$;
\end{itemize}
\item[2)] $A_{14|3}^{k}$ with $\Delta_{14|3}^{k}(x)=1 \otimes x+ x\otimes 1- x\otimes x+z\otimes z$ for $k=5,6$:
\begin{itemize}\itemsep=0pt
\item $\Delta_{14|3}^{5}(y)=1\otimes y+y\otimes 1$,
\item $\Delta_{14|3}^{6}(y)=1 \otimes y+ y\otimes 1+
x\otimes z+ z\otimes x$;
\end{itemize}
\item[3)] $A_{14|3}^{7}=(A_{14|3}^{3})^{\rm cop}$.
\end{enumerate}

\textit{Superalgebra $A_{15|3}$}, we have
$A_{15|3}^{1}=(A_{14|3}^{1})^{\rm op}$, $A_{15|3}^{2}=(A_{14|3}^{2})^{\rm op}$,
$A_{15|3}^{3}=(A_{14|3}^{3})^{\rm op}$, $A_{15|3}^{4}=(A_{14|3}^{4})^{\rm op}$,
$A_{15|3}^{5}=(A_{14|3}^{5})^{\rm op}$, $A_{15|3}^{6}=(A_{14|3}^{6})^{\rm op}$,
$A_{15|3}^{7}=(A_{14|3}^{5})^{\rm op,cop}$.

\textit{Superalgebra $A_{17|2}$}, we have:
\begin{enumerate}\itemsep=0pt
\item[1)] $A_{17|2}^{1}$ with $\Delta_{17|2}^{1}(x)=1\otimes x+x\otimes 1-x\otimes x$,
$\Delta_{17|2}^{1}(y)=1\otimes y+y\otimes 1-y\otimes x-x\otimes y$,
\\
$\Delta_{17|2}^{1}(z)=1\otimes z+z\otimes 1-z\otimes x-x\otimes z$, $\varepsilon_{15|3}^{7}(x)=0$;

\item[2)] $A_{17|2}^{2}$ with $\Delta_{17|2}^{2}(x)=x\otimes x$, $\Delta_{17|2}^{2}(y)=x \otimes y+y\otimes x$,
$\Delta_{17|2}^{2}(z)=x\otimes z+ z\otimes x$,
$\varepsilon_{15|3}^{7}(x)=1$.
\end{enumerate}

\subsection*{Acknowledgements}

The authors are grateful to the referees for their valuable remarks and suggestions.

\pdfbookmark[1]{References}{ref}
\LastPageEnding

\end{document}